%% file: TomV19.tex
\documentclass[a4paper,10pt,notitlepage]{article}
\usepackage{fullpage}
\usepackage{amsthm}
\usepackage{amsmath,amssymb}
\usepackage{changepage}
\usepackage{enumerate}
\usepackage{amsfonts}
\usepackage{graphics,graphicx}
\usepackage{wrapfig}
\usepackage{wrapfig,multirow,rotating}

\usepackage{pgfplots,tikz,subfigure}
\usepackage{amsmath,amsthm,amssymb}
\usepackage{hyperref,graphics,graphicx,color,algorithm,algorithmic}
\usepackage{mymacros,wrapfig,multirow,rotating}
\newcommand{\ri}{{\mathrm{i}}}

\DeclareMathOperator{\hin}{\mathcal{H}_\infty}
\DeclareMathOperator*{\argminx}{arg\,min}
\newcommand{\T}{\mathsf{T}}
\renewcommand{\H}{\mathsf{H}}
\renewcommand{\hat}{\widehat}
\renewcommand{\tilde}{\widetilde}
\title{Semi-active $\cH_{\infty}$ damping optimization by adaptive interpolation\thanks{This work was supported by a travel grant of the \emph{Ministry of Science and Education of the Republic of Croatia (MZOS)} and the \emph{German Academic Exchange Service (DAAD)} and by Croatian Science Foundation under the projects, Grant No. IP-2014-09-9540 and IP-2019-04-6774.
}}
\author{{Zoran Tomljanovi\'{c}\thanks{Department of Mathematics, Josip Juraj Strossmayer University of Osijek, Trg Ljudevita Gaja 6, Osijek, HR-31000, Croatia (e-mail: ztomljan@mathos.hr).}} and Matthias Voigt\thanks{University of Hamburg, Department of Mathematics, Center for Optimization and Approximation, Bundesstra{\ss}e 55, D-20146 Hamburg, Germany (e-mail: matthias.voigt@uni-hamburg.de) and Technische Universit\"at Berlin, Institut f\"ur Mathematik, Stra{\ss}e des 17. Juni 136, D-10623 Berlin, Germany (e-mail: mvoigt@math.tu-berlin.de).}}

\theoremstyle{plain}
\newtheorem{mydef}{Definition}%[section]

\newtheorem{lemma}[mydef]{Lemma}
\newtheorem{ass}[mydef]{Assumption}
\newtheorem{example}{Example}

\theoremstyle{definition}
\usepackage{algorithm,algorithmic}

\DeclareMathOperator{\diag}{diag}

\begin{document}
\maketitle

\begin{abstract}
 In this work we consider the problem of semi-active damping optimization of mechanical systems with fixed damper positions. Our goal is to compute a damping that is locally optimal with respect to the $\hin$-norm of the transfer function from the exogenous inputs to the performance outputs. We make use of a new greedy method for computing the $\hin$-norm of a transfer function based on rational interpolation. In this paper, this approach is adapted to parameter-dependent transfer functions. The interpolation leads to parametric reduced-order models that can be optimized more efficiently. At the optimizers we then take new interpolation points to refine the reduced-order model and to obtain updated optimizers. In our numerical examples we show that this approach normally converges fast and thus can highly accelerate the optimization procedure. Another contribution of this work are heuristics for choosing initial interpolation points.
\end{abstract}

\section{Introduction}
Consider a vibrational system described by a system of a second-order differential equations
\begin{subequations}\label{eq:MDKs}
\begin{align}
  M\ddot{q}(t)+C_{\rm int} \dot q(t)+Kq(t) &= B_2 u(t)+E_2w(t), \label{MDK1s} \\
                              y(t) &= C_2\dot q(t), \label{MDK2s} \\
                              z(t) &= H_1q(t), \label{MDK3s}
\end{align}
\end{subequations}
where $M,\,K \in \R^{n \times n}$ are the symmetric and positive definite mass and stiffness matrices, respectively.  The internal damping $C_{\rm int} \in \R^{n \times n}$ is usually taken to be a small multiple of the critical damping denoted by
$C_{\rm crit}$, that is
\begin{equation}\label{Cint}
C_{\rm int} = \alpha_c C_{\rm crit},
\end{equation}
see, e.\,g., \cite{BennerTomljTruh10,BennerTomljTruh11,TRUHVES09}. The presence of the positive definite internal damping $C_{\rm int}$ ensures that the homogeneous system \eqref{MDK1s} is asymptotically stable.
%\begin{equation}\label{C_{\rm int}}
%  C_{\rm int} = \alpha_c C_{crit},\quad \mbox{where} \quad C_{crit} =2 M^{1/2}\sqrt{M^{-1/2}KM^{-1/2}}M^{1/2}.
%\end{equation}
Moreover, we have $B_2\in \mathbb{R}^{n\times p}$, $E_2\in \mathbb{R}^{n\times m}$, $C_2\in \mathbb{R}^{p\times n}$, and $H_1 \in \mathbb{R}^{\ell \times n}$.

The vector $q(t) \in \mathbb{R}^n$ represents the displacements of the masses, $u(t) \in \R^p$ is the control input, and $y(t) \in \R^p$ is the measured output. Moreover, $w(t) \in \Rm$ is an exogenous input signal, while $z(t) \in \mathbb{R}^\ell$ denotes the performance output.
%If one is  also interested in the velocities, the
In principal, $z(t)$ can also include an additional part which corresponds to the velocities of the masses $\dot{q}(t)$, but in this paper only the states are of interest.

%The vector   $w(t)\in \mathbb{R}^q$ is the   primary excitation (i.e. noise) input.
%On the other hand the vector $u(t)\in \mathbb{R}^p$ determines control and we   consider the particular case of a
In this paper we assume negative linear  feedback corresponding to a linear damper of the form
\begin{equation}\label{damper form}
u(t) = -Gy(t),
\end{equation}
where $G = \diag(g_1, g_2, \ldots, g_p)\in \mathbb{R}^{p\times p}$ is the
diagonal damping matrix with non-negative parameters encoded in a vector $g =
\begin{pmatrix} g_1 & g_2 & \ldots & g_p\end{pmatrix}^\T$. The entries $g_i$,
$i=1,\,2,\,\ldots,\,p$ represent the friction coefficients of the corresponding
dampers, usually called gains. In this paper we will consider the gains to be
scalar variables which need to be optimized.

Here we assume that we have collocated inputs and outputs, i.\,e., $C_2 = B_2^\T$. By using the feedback control law as in \eqref{damper form} we obtain the \emph{closed-loop system} (For simplicity, we omit the dependency of $q(\cdot)$ on $g$ and only add it to $z(\cdot)$.)
\begin{subequations}\label{eq:MDK}
\begin{align}
  M\ddot q(t)+C(g)\dot q(t)+Kq(t) &= E_2 w(t)\label{MDK1},\\
                        z(g,t) &= H_1 q(t)\label{MDK3},
\end{align}
\end{subequations}
where
\begin{equation*}%\label{Cg}
C(g):=  C_{\rm int} + B_2\diag(g_1, g_2, \ldots, g_p)B_2^\T
\end{equation*}
is symmetric and positive definite. This implies that the unforced closed-loop system is also asymptotically stable. More details regarding system stability and model description can be found in \cite{BKTT15, Blanchini12, BennerTomljTruh11}.

With the substitutions $x_1(t):=q(t)$, $x_2(t):=\dot q(t)$ and $x(t):=\left[%
\begin{smallmatrix}
  x_1(t) \\
  x_2(t)
\end{smallmatrix}%
\right]$ we obtain a first-order representation of the closed-loop system
\begin{align*}
\begin{bmatrix} I_n & 0 \\ 0 & M \end{bmatrix} \dot x(t) &=
\begin{bmatrix}
  0 & I_n \\
  - K & - C(g)
\end{bmatrix} x(t) +
\begin{bmatrix}
  0 \\
   E_2
\end{bmatrix} w(t), \\
z(g,t) &= \begin{bmatrix} H_1 & 0 \end{bmatrix} x(t).\nonumber
\end{align*}
%where
%\begin{align}\label{matrix A}
%x=\left[%
%\begin{array}{c}
%  x_1 \\
%  x_2 \\
%\end{array}%
%\right], \quad  A(g)&=\left[%
%\begin{array}{cc}
%  0 & I \\
%  - K & - C(g) \\
%\end{array}%
%\right],\quad
%E=\left[%
%\begin{array}{c}
%  0 \\
%   E_2 \\
%\end{array}%
%\right],\\
% \mbox{and} \quad  H&=\left[
%         \begin{array}{cc}
%           H_1 & 0\\
%         \end{array}
%       \right].\label{matrix H}
%\end{align}

%The matrices $B_2$ and $G$ constitute the semi-active damping part of the  above
%second order system.

%\eqref{MDK1s}-\eqref{MDK3s} is asymptotically stable which in  the usual assumption. For more details, see, e.g. More details regarding system stability and model description can be found in \cite{BKTT15, Blanchini12, BennerTomljTruh11}.

Using the Laplace transform we obtain the closed-loop transfer function of \eqref{eq:MDK}, which is given by
\begin{align}\label{TF}
 \begin{split}
   F(g,s) &= H_1 \left(s^2M +sC(g) + K\right)^{-1} E_2 \\
          &=  \begin{bmatrix}
           H_1 & 0
         \end{bmatrix}
            \left(s\begin{bmatrix} I_n & 0 \\ 0 & M \end{bmatrix} -
       \begin{bmatrix}
  0 & I_n \\
  - K & - C(g)
\end{bmatrix}\right)^{-1} \begin{bmatrix}
  0 \\
   E_2
\end{bmatrix} =: \cC \cD(g,s)^{-1} \cB,
\end{split}
\end{align}
which is a real-rational matrix-valued function for each $g \in (\R_{\ge 0})^p$, where $\R_{\ge0}$ denotes the set of nonnegative real numbers.
%where in transfer function we have also parameter $g$ which is parameter of interest for the damping optimization setting.

The damping optimization problem has been investigated widely in the literature, for example in the books \cite{Gaw, VES2011, MullerSchiehlen85} where one can also find overview of different damping criteria. But even nowadays there are lot of open problems, especially from the computational point of view.
One case studied in the past is the \emph{passive damping}. This means that for given positive definite matrices $M$ and $K$, the ``best'' damping matrix $C(g)$ (with a certain structure) should be determined, such that the solution trajectories of the stationary system
\begin{equation*}
 M\ddot q(t) + C(g) \dot q(t) + Kq(t) = 0
\end{equation*}
have an ``optimal'' transient behavior. In this setting, gain optimization can be a computationally very demanding task, especially if the damper parameter values should be optimized for many different damper locations. There exist a number of different methods such as \cite{FreitLanc:99,BRAB98}. Moreover, there are methods based on dimension reduction such as  \cite{BennerTomljTruh10, BennerTomljTruh11, TrTomPuv16, TRUTOMVES2014}, where the minimization of the total average energy is the objective.
%A review of different methods for this dimension reduction can be found in \cite{ANT05, BenMS05, BennerGW15}.
%Passive damping setting was studied in   \cite{BennerTomljTruh10, BennerTomljTruh11, TrTomPuv16, TRUTOMVES2014, BRAB98} where authors considered optimization using criterion of total average energy, while criteria that involve spectral abscissa were considered  e.g. in  \cite{FreitLanc:99}.

%%%
In this work we consider the \emph{semi-active damping optimization problem:} we want \textit{to determine ``the best'' damping matrix $C(g)$ which will minimize the output $z(g,\cdot)$ under the influence of the input $w(\cdot)$}.
The influence of the input on the output can be measured in several ways. There exist several optimization criteria such as the impulse response energy, the peak-to-peak performance, and the energy-to-energy performance. An overview of the different criteria can be found in \cite{Blanchini12}.

%which are often based on system norms (see, e.\,g., \cite{Blanchini12}).

%One criterion is the impulse response
%energy which can be shown to correspond to the $\mathcal{H}_2$-norm of the system's transfer function, see  \cite{Blanchini12, Burl1998, Zhou96}.

%Energy-to-energy gain

% In the frequency domain we can
% write impulse response  cost function $J_2$ in terms of  the  transfer function  matrix
% \eqref{TF} via
% \begin{equation}\label{crit in terms TF}
%  \left\| F(g,\cdot) \right\|_{\mathcal{H}_2} =\int_{-\infty}^{+\infty} \trace{ F(g,j\omega)^*F(g,j\omega) }
% \mathrm{d}\omega.
% \end{equation}
% Moreover impulse response energy can be characterized in terms of the solution of the Lyapunov equation where more details can be found in  e.g., \cite{Blanchini12, Burl1998, Zhou96}.

Damping optimization using the impulse response energy leads to the minimization of the $\cH_2$-norm of the closed-loop transfer function. This requires to solve an associated Lyapunov equation numerous times, which can make the optimization process computationally inefficient. In order to accelerate it, in \cite{BKTT15} the parametric (subspace accelerated) dominant pole algorithm is used for the approximation of the impulse response energy. Moreover, in \cite{TomljGugerBeatt17} the authors propose an efficient optimization approach using structure-preserving parametric model reduction based on the iterative rational Krylov algorithm (IRKA). There, several adaptive sampling strategies are used to obtain good approximations with respect to the $\cH_2$-norm, aligning well with the underlying design objectives.

%Since this approach depends on parameter sampling that occurs within each parametric model reduction cycle,

%authors considered an approaches with predetermined sampling and approaches using adaptive sampling.

The method used in this paper is also based on adaptive sampling, but it differs from the approach for the $\cH_2$-norm case. However, here we consider the energy-to-energy performance of the closed-loop system which is defined by
\begin{equation*}
 J_\infty(g) = \sup_{\substack{ w \in {\mathcal{L}_2(\R_{\ge 0},\R^\ell) \setminus \{ 0 \}}, \\  x(0) = 0}} \frac{\left\| z(g,\cdot) \right\|_{\mathcal{L}_2(\R_{\ge 0},\R^m)}}{\left\| w \right\|_{\mathcal{L}_2(\R_{\ge 0},\R^\ell)}},
\end{equation*}
where $(w,\,x,\,z(g,\cdot))$ are measurable and square-integrable solution trajectories of the closed-loop system with parameters $g$.
Hence, the energy-to-energy performance can be interpreted as the worst-case amplification of the energy of an exogenous input signal in the performance output. It is well-known (see, e.\,g., \cite[Chap. 3]{Zhou96}) that this criterion is equivalent to the $\cH_\infty$-norm of the closed-loop transfer function, that is
\begin{equation*}
 J_\infty(g) = \left\| F(g,\cdot) \right\|_{\cH_\infty}.
\end{equation*}
%On the other hand in this paper we would like to consider another important optimization criterion which is based on $\mathcal{H}_\infty$ norm of the considered system. This criterion is a more challenging from damping optimization setting, but it also gives some advantages in the terms of optimal gains as we will illustrate in this paper.
Since the closed-loop system is asymptotically stable, then for all $g \in (\R_{\ge 0})^p$, the functions $F(g,\cdot)$ are in the space
\begin{equation*}
  \cH_\infty^{m \times \ell} := \left\{ F : \C^+ \rightarrow {\mathbb C}^{m \times \ell} \; \bigg| \; F \text{ is analytic and } \sup_{\lambda \in \C^+} \left\| F(\lambda) \right\|_2 < \infty \right\},
\end{equation*}
where $\C^+ := \{ \lambda \in \C \;|\; \Real{\lambda} > 0 \}$. This space is equipped with the $\cH_\infty$-norm
\begin{equation*}
 \left\| F \right\|_{\cH_\infty} := \sup_{\lambda \in \C^+} \left\| F(\lambda)
\right\|_2 = \sup_{\omega \in \R} \left\| F(\ri\omega) \right\|_2 = \sup_{\omega
\in \R} \sigma_{\max}(F(\ri \omega)),
\end{equation*}
see \cite[Chap. 3]{Zhou96}.

% In this paper we will some
%
% In order to write our cost function in terms of $\mathcal{H}_\infty$ system norm, first we note that transfer function  can be written as
% \begin{equation*}
%  F(g,s) = \left[
%          \begin{array}{cc}
%            H_1 & 0\\
%          \end{array}
%        \right]\left(s\begin{bmatrix} I_n & 0 \\ 0 & M \end{bmatrix}-\left[%
% \begin{array}{cc}
%   0 & I \\
%   - K & - C(g) \\
% \end{array}%
% \right]\right)^{-1} \left[%
% \begin{array}{c}
%   0 \\
%    E_2 \\
% \end{array}%
% \right],
% \end{equation*}
% where $C(g)$ is given by \eqref{Cg}.
%
%
% In compact form we can write transfer function as
% \begin{equation}
% \label{TF1}
% F(g,s) =C(sE-A(g))^{-1} B,
% \end{equation}
% with
% \begin{align*}
%  C&=\vphantom{\begin{bmatrix} 0 	\\ 	E_2 \end{bmatrix}}\begin{bmatrix} H_1 & 0 \end{bmatrix},\quad  B =\begin{bmatrix} 0 	\\ 	E_2 \end{bmatrix},\\
%  E&= \begin{bmatrix} I_n	& 	0	\\	0 	& 	M \end{bmatrix},\quad  A(g) = \begin{bmatrix} 0	&	I_n	\\	-K 	& 	-C(g) \end{bmatrix}.
% \end{align*}
%
% Then $\mathcal{H}_\infty$  norm of a system can be written as
% \[
%  \left\| F(g,\cdot) \right\|_{\cH_\infty} =\sup_{\omega\geq 0} |F(g,j\omega)|.
% \]
% Moreover it can be shown that
% \begin{equation*}
% \displaystyle \left\| F(g,\cdot) \right\|_{\cH_\infty} = \sup_{\omega\geq 0}  \sigma_{\max}\left( F(g,\ri\omega ) \right).
% \end{equation*}
Thus, in the setting of $\cH_\infty$ damping optimization the problem is to determine optimal gains
\begin{equation*}
  g_* = \argminx_{g \in (\R_{\ge 0})^p} \left\| F(g,\cdot) \right\|_{\cH_\infty}
\end{equation*}
(under the assumption that the minimum exists).

The optimization of damping based on the closed-loop $\mathcal{H}_\infty$-norm can be a computationally challenging task.
First, the optimization problem we consider here is nonlinear, nonconvex, and nonsmooth. The latter means that $g_*$ may be attained at a point, where the objective function $g \mapsto \left\| F(g,\cdot) \right\|_{\cH_\infty}$ is not differentiable. This problem occurs frequently in fixed-order $\cH_\infty$-control, see \cite{GumHMO09}. Specialized methods for such optimization problems are available \cite{LewO08,LewO13,CurtisMitchOverton17}. These are modifications of quasi-Newton methods but in case of a nonsmooth optimizer, the $\cH_\infty$-norm (and its gradient) may have to be evaluated for a lot of parameter values.

Moreover, the calculation of the $\mathcal{H}_\infty$-norm for a fixed parameter can be expensive, especially if the state-space dimension of the closed-loop system is large. The latter problem has been addressed by various works \cite{GugGO13,BenV13e,MitO16}. However, in this paper we will modify the idea from \cite{AliBMSV17} to our problem. In the latter work, the transfer function is approximated by interpolatory reduced-order models. Then the projection spaces are updated using the information of the point on the imaginary axis, at which the
maximum of the largest singular value for the reduced order model is attained. This leads to an iterative algorithm that can be shown to have a superlinear rate of convergence to a local maximizer. In \cite{AliBMV18}, this approach has been extended to the problem of minimizing the $\cH_\infty$-norm of a parameter-dependent transfer function which is particularly useful in our context. With this approach, in each step we construct a parametric reduced-order model with reduced transfer function $\tilde{F}(g,s)$ and compute $\hat{g} \in (\R_{\ge 0})^p$ and $\hat{\omega} \in \R \cup\{-\infty,\infty\}$ such that
\begin{align}\label{eq:optimal}
 \begin{split}
 \big\|\tilde{F}\big(\hat{g}, \cdot \big)\big\|_{\cL_\infty} &= \min_{g \in (\R_{\ge 0})^p} \big\| \tilde{F}(g, \cdot) \big\|_{\cL_\infty}, \\
 \big\|F\big(\hat{g},\ri\hat{\omega}\big)\big\|_2 &= \max_{\omega \in \R \cup\{-\infty,\infty\}} \big\| F\big(\hat{g},\ri \omega\big) \big\|_2.
 \end{split}
\end{align}
Hence, in each step, one optimization   of the $\cL_\infty$-norm (for the definition see \eqref{LinfDef}) of a parametric reduced transfer function and one large-scale non-parametric $\cH_\infty$-norm computation is needed.
Then we choose $\big(\hat{g},\ri\hat{\omega}\big)$ as a new interpolation point to update the reduced-order model and repeat this process until convergence.

%The third issue in computational efficiency is that the problem is to optimize damping positions in general. Efficient optimization of damping positions is a very demanding problem and standard approaches  require  optimization of gains for a large number of different damping positions. This means that even for moderate dimensions, damping optimization using $\mathcal{H}_\infty$-norm will be computationally very demanding since that we need to calculate $\mathcal{H}_\infty$-norm for a large number of different parameters. Thus, in this paper we would like to accelerate this optimization process by  adaptive interpolation.

%Therefore, in this paper we will propose an approach which will provide efficient approximation of $\mathcal{H}_\infty$ norm which will be efficiently used for   acceleration of damping optimization.

However, in the context of damping optimization, the choice of the interpolation points is particularly subtle. Therefore, the main contribution of this work consists of deriving two heuristics for choosing initial interpolation points. We further do not only demonstrate the effectiveness, but also the numerical challenges of these ideas on various numerical examples.

%In this paper, first we would like to present a formula for a gradient of a objective function $F(g,\cdot)$. Moreover our objective function may not be differentiable in some points, but still we can optimize gains  using approach that is based on Quasi-Newton minimization procedure which can significantly accelerate gain optimization. Furthermore, we would like to compare criterion of $ \left\| F(g,\cdot) \right\|_{\cH_\infty}$ norm with other relevant criteria.

\section{Preliminaries}\label{sec:eval}
In this section we discuss some preliminary results for our method. In
particular, we give some more details on how to compute the $\cH_\infty$-norm
since this is the core idea of our algorithm. This method is derived in detail
in \cite{AliBMSV17} and we closely follow the presentation in
\cite{AliBMSV17a}.
Moreover, we derive an analytic formula for the gradient of the
$\cH_\infty$-norm with respect to the parameters $g$.
\subsection{Computation of the $\cH_\infty$-Norm}\label{subsec:hinf}
Assume that $F \in \cH_\infty^{m \times \ell}$ with
\begin{equation*}
 F(s) = C(sE-A)^{-1}B
\end{equation*}
is given. To compute $\left\| F \right\|_{\cH_\infty}$, in \cite{AliBMSV17}
reduced functions
\begin{equation}\label{eq:transfer_func_red}
      \widetilde{F}(s) := \widetilde{C} \big(s\tilde{E} - \tilde{A}\big)^{-1} \widetilde{B},
\end{equation}
with
\begin{equation}\label{eq:middle_factor}
 s \tilde{E} - \tilde{A} = s W^\H E V - W^\H A V, \quad \tilde{B} = W^\H B, \quad \tilde{C} = CV,
\end{equation}
and projection matrices $V,\,W \in \C^{n \times k}$ where $k \ll n$ are
iteratively constructed. Since $s \tilde{E} - \tilde{A}$ is of very small
dimension compared to $sE-A$, the $\cL_\infty$-norm of $\tilde{F}$ can be
easily obtained. Assume that $F$ is a transfer function that is bounded on the
imaginary axis. Then the $\mathcal{L}_\infty$-norm of $F$ is defined by
\begin{equation}\label{LinfDef}
 \left\| F \right\|_{\cL_\infty} = \sup_{\omega \in \R} \left\| F(\ri\omega) \right\|_2 = \sup_{\omega \in \R} \sigma_{\max}(F(\ri \omega)).
\end{equation}
Note that, if $F \in \cH_\infty^{m \times \ell}$, then we have $\left\| F \right\|_{\cH_\infty} = \left\| F \right\|_{\cL_\infty}$.

In the method a sequence of such reduced functions
$\widetilde{F}_1,\,\widetilde{F}_2,\,\ldots$ is constructed, whose
$\cL_\infty$-norms converge to the $\cH_\infty$-norm of the original
transfer function $F$. To obtain the reduced functions, the
interpolation technique of \cite{Gugercin2009} is employed. If
$m=\ell$ and $r$ interpolation points $\ri \omega_1,\,\ldots,\,\ri\omega_r$
are given, then $\widetilde{F}_r$ is defined by $\widetilde{F}$
as in \eqref{eq:transfer_func_red} and \eqref{eq:middle_factor} with $V = V_r$, $W=W_r$ set to
\begin{equation}\label{subspaces}
 V_r = \begin{bmatrix} (\ri\omega_1 E - A )^{-1}B & \ldots & (\ri\omega_r E - A )^{-1}B \end{bmatrix}, \quad
 W_r = \begin{bmatrix} (\ri\omega_1 E - A )^{-\H}C^\H & \ldots & (\ri\omega_r E - A)^{-\H}C^\H \end{bmatrix}.
\end{equation}
These choices for $V_r$ and $W_r$ imply that the Hermite interpolation properties
\begin{equation}\label{Hermite}
 F(\ri\omega_k) = \widetilde{F}_r(\ri\omega_k), \quad F'(\ri\omega_k) = \widetilde{F}'_r(\ri\omega_k), \quad k = 1,\,\ldots,\,r
\end{equation}
are satisfied.
It is important that  $V_r$ and $W_r$ have the same number of columns so that the inverse 
$\big(s \tilde{E} -\tilde{A} \big)^{-1}$ exists. This condition is usually
violated by the projection matrices in \eqref{subspaces} if $m \neq \ell$. By using the projection matrices
\begin{align*}
 V_r &= \begin{bmatrix} (\ri\omega_1 E - A )^{-1}B F(\ri \omega_1)^\H & \ldots & (\ri\omega_r E - A )^{-1}B F(\ri \omega_r)^\H \end{bmatrix},
 \quad W_r \text{ as in \eqref{subspaces}}, \quad
 \text{ if } \ell > m, \\
 W_r &= \begin{bmatrix} (\ri\omega_1 E - A)^{-\H}C^\H F(\ri \omega_1) & \ldots & (\ri\omega_r E - A )^{-\H} C^\H F(\ri \omega_r) \end{bmatrix},
 \quad V_r \text{ as in \eqref{subspaces}}, \quad
 \text{ if } \ell < m
\end{align*}
instead, this problem is solved and the Hermite
interpolation properties \eqref{Hermite} are preserved, see \cite{AliBMSV17} for the
derivation. Alternatively, a regularization procedure may be carried out to
achieve this goal.
Now, since the transfer function $\widetilde{F}_r$ is constructed by matrices of small dimension, its $\mathcal{L}_\infty$-norm can be obtained by well-established methods for the small and dense case which are reported in \cite{BS90,BB90,BenSV12}. With the
point $\ri \omega_{r+1}$ at which the $\cL_\infty$-norm of the current iterate
$\tilde{F}_r$ is attained, the projection matrices are updated as
\begin{align*}
 V_{r+1} &:= \begin{bmatrix} V_r & (\ri\omega_{r+1} E - A )^{-1}B \end{bmatrix} \quad (\text{or } V_{r+1} := \begin{bmatrix} V_r & (\ri\omega_{r+1} E - A)^{-1}B F(\ri \omega_{r+1})^\H \end{bmatrix}), \\
 W_{r+1} &:= \begin{bmatrix} W_r & (\ri\omega_{r+1} E - A)^{-\H}C^\H \end{bmatrix} \quad (\text{or } W_{r+1} := \begin{bmatrix} W_r & (\ri\omega_{r+1} E - A )^{-\H}C^\H F(\ri\omega_{r+1}) \end{bmatrix}).
\end{align*}
In case of convergence, using the Hermite interpolation conditions \eqref{Hermite}, a superlinear convergence rate to a local maximizer of $\sigma_{\max}(F(\ri\cdot))$ can be shown, see \cite{AliBMSV17} for the details. A MATLAB implementation of this procedure, called \texttt{linorm\_subsp}, is publicly
available\footnote{downloadable from
\url{http://www.math.tu-berlin.de/index.php?id=186267&L=1}.}.

% In this section we will propose a projection method to obtain approximation of the $\cH_\infty$ norm of transfer function $F(g,s)$  given by \eqref{TF1}. Our purpose is to compute two subspaces $\cU$ and $\cV$ formally represented by span of the columns of the two matrices $U$ and $V$, respectively. In order to have efficient approximation, these subspaces  will need to have dimensions much smaller than $n$ and to construct the reduced order model, we enforce the Petrov-Galerkin condition such as the reduced state at time $t$ is given by $U\tx(t)$. Then, the reduced system is defined by
% \bea
% V^*(EU\dot{\tx}(t) - A(g)U\tx(t) - Bu(t)) = 0 \quad \text{and} \quad \ty(t) = CU\tx(t).
% \eea
% The transfer function associated with the reduced  system is
% \bea
% \label{reduce}
% \tF(g,s) = \underbrace{CU}_{=:\tC}{\underbrace{(sV^*EU-V^*A(g)U)}_{=:(s\tE-\tA(g))^{-1}}}^{-1}\underbrace{V^*B}_{=:\tB}.
% \eea
%
% In our approach we will construct the matrices $U$ and $V$ iteratively using the following steps:
% \begin{enumerate}
% 	\item In each iteration and for a fixed damping parameter $g$, compute the optimal frequency for which the $\mathcal{H}_\infty$ of the transfer function is obtained.
% 	\item Use this optimal frequency to expand the size of the subspaces $U$ and $V$ as in \cite{AliBMSV17}.
% 	\item Reduce the transfer function with $U$ and $V$.
% 	\item Compute an optimal damping parameter $\hg$ with respect to the reduced transfer function.
% 	\item If there is no convergence, then use this parameter in the first step and continue.
% \end{enumerate}

\subsection{Gradient of the $\cH_\infty$-Norm}
For the derivation of the gradient of the $\cH_\infty$-norm we make use of the following lemmas:
\begin{lemma} \cite{Lan64}
 Let $L : \R \to \C^{m \times \ell}$ with $L(t) = L_0 + tL_1$ for some
$L_0,\,L_1 \in \C^{m\times \ell}$. Let $\sigma(t)$ be a singular value of
$L(t)$ converging to a simple nonzero singular value $\sigma_0$ of $L_0$ as $t
\rightarrow 0$. Then, $\sigma(t)$ is analytic near $t=0$ and
 \begin{equation*}
  \left.\frac{\mathrm{d}\sigma(t)}{\mathrm{d}t}\right|_{t=0} =
\operatorname{Re}\big(v_0^\H L_1 u_0\big),
 \end{equation*}
 where $u_0$ and $v_0$ with $\left\|u_0\right\|_2 = \left\|v_0\right\|_2 = 1$ are, respectively, the right and left singular vectors of $L_0$ corresponding to $\sigma_0$.
\label{lem:svpert}
\end{lemma}
\begin{lemma} \cite{MS05}
 Let $s_0 \in \C$ not be a pole of the transfer function $F(s)=C(sE-A)^{-1}B$. Then, $F(\cdot)$ can be expanded into a Laurent series at $s_0$ as
 \begin{align*}
  F(s) &= C(I_n - \left(s-s_0)(s_0E-A)^{-1}E\right)^{-1}(s_0E-A)^{-1} B \\
       &= C(s_0E-A)^{-1}B - C(s_0E-A)^{-1}E(s_0E-A)^{-1}B (s-s_0) + \mathcal{O}\left((s-s_0)^2\right).
 \end{align*}
 \label{lem:laurent}
\end{lemma}
From the two lemmas above we obtain the following result.
\begin{lemma} \cite{BenV13e}
 Let $s_0 = \gamma_0 + \ri \delta_0 \in \C$ not be a pole of the transfer
function $F(\cdot)$. Furthermore, assume that the largest singular value of
$F(s_0)$ is simple with associated right and left singular vectors $u_0$ and
$v_0$ satisfying $\left\|u_0\right\|_2 = \left\|u_0\right\|_2 = 1$ and define
$h(\gamma,\delta) := \left\|F(\gamma+ \ri \delta)\right\|_2$. Then the gradient
of $h$ with respect to the variable $\gamma$ is given by
 \begin{align} \label{eq:gradform}
 \nabla_\gamma h(\gamma_0,\delta_0) =
-\operatorname{Re}\left(v_0^\H
C\left(s_0E-A\right)^{-1}E\left(s_0E-A\right)^{-1}Bu_0\right).
 \end{align}
 \label{lem:nabla}
\end{lemma}
Now we derive a gradient formula for the function $g \mapsto \left\|
F(g,\cdot) \right\|_{\cH_\infty}$. For the derivation we make the following
assumptions.
% Thus, let us recall that first from \eqref{TF1} for transfer function it holds that
% \begin{equation*}
%  F(g,s) = \begin{bmatrix} H_1 & 0 \end{bmatrix} \left(s \begin{bmatrix} I_n & 0 \\ 0 & M \end{bmatrix} - \begin{bmatrix} 0 & I_n \\ - K & -C_{\rm int} - B_2 G B_2^T \end{bmatrix} \right)^{-1} \begin{bmatrix} 0 \\ E_2 \end{bmatrix},
% \end{equation*}
% with $g = \begin{bmatrix} g_1 & \ldots & g_m \end{bmatrix}^T$ and $G = \diag {g_1,\, \ldots,\, g_m}$.
\begin{ass} \label{ass:ass}
Assume that for some fixed $g_0$, we have
\begin{equation*}
 \left\| F(g_0,\cdot) \right\|_{\cH_\infty} = \left\| F(g_0,\ri\omega_0) \right\|_2 = \sigma_{\max}\left( F(g_0,\ri\omega_0) \right)
\end{equation*}
for some $\omega_0 \in \R_{\ge 0} \cup\{\infty\}$. Assume further that:
 \begin{enumerate}
  \item[A1.] The point $\omega_0 \in \R_{\ge 0} \cup \{\infty\}$ at which the $\cH_\infty$-norm is attained, is unique.
  \item[A2.] The maximum singular value of $F(g_0,\ri\omega_0)$ is simple with associated right and left singular vectors $u_0$ and $v_0$ satisfying $\left\|u_0\right\|_2 = \left\|v_0\right\|_2 = 1$.
 \end{enumerate}
\end{ass}
 Assumption A1 ensures that the optimal frequency $\omega_0$ at which the $\cH_\infty$-norm is attained is also uniquely determined in a neighborhood of $g_0$. Therefore, in a neighborhood of $g_0$ there are no jumps in the points at which the $\cH_\infty$-norm is attained and thus there are no ``kinks'' in the graph of the function $g \mapsto \left\| F(g,\cdot) \right\|_{\cH_\infty}$. Assumption A2 ensures that we can apply Lemma~\ref{lem:nabla} to compute the gradient of the maximum singular value function using the corresponding singular vectors. Both conditions together imply that the function $g \mapsto \left\| F(g,\cdot) \right\|_{\cH_\infty}$ is differentiable in a neighborhood of $g_0$ and that we can make of use of formula \eqref{eq:gradform} to compute the gradient of the $\cH_\infty$-norm as outlined below.

For $g = \begin{pmatrix} g_1 & g_2 & \ldots & g_p\end{pmatrix}^\T$ and all $j = 1,\,\ldots,\,p$, we define the matrices $$G_j = \operatorname{diag}(g_1,\, \ldots,\, g_{j-1},\,0,\,g_{j+1},\,\ldots,\, g_p) \quad \text{and} \quad F_j = \frac{1}{g_j}(G-G_j).$$
Then, we can rewrite $F(g,s)$ as
\begin{align*}
 F(g,s) &= \begin{bmatrix} H_1 & 0 \end{bmatrix} \left(g_j \begin{bmatrix} 0 & 0 \\ 0 & B_2F_jB_2^\T \end{bmatrix} - \begin{bmatrix} -sI_n & I_n \\ -K & -sM-C_{\rm int} - B_2 G_j B_2^T \end{bmatrix} \right)^{-1} \begin{bmatrix} 0 \\ E_2 \end{bmatrix} \\
        &=: \cC(g_j\cE_j - \cA_j(s))^{-1}\cB, \qquad j = 1,\, \ldots,\, p.
\end{align*}
Now, let
\begin{equation*}
 \left\| F(g,\cdot) \right\|_{\cH_\infty} = \sigma_{\max}\left( F(g,\ri\omega_0) \right)
\end{equation*}
for some $\omega_0 \in \R_{\ge 0}$ and let the corresponding right and left singular vectors be $u_0 \in \C^{\ell}$ and $v_0 \in \C^m$. Let Assumption~\ref{ass:ass} hold accordingly for $g_0 = g$.
Then with Lemma \ref{lem:nabla} we obtain
\begin{align}\nonumber
 \nabla_g \left\| F(g,\cdot) \right\|_{\cH_\infty} &= \begin{bmatrix} {\partial \sigma_{\max}\left(F(g,\ri\omega_0)\right)}/{\partial g_1} \\ \vdots \\ {\partial \sigma_{\max}\left(F(g,\ri\omega_0)\right)}/{\partial g_p}\end{bmatrix} \\ &= \begin{bmatrix} -\operatorname{Re}\left(v_0^\H \cC\left(g_1\cE_1-\cA_1(\ri\omega_0)\right)^{-1}\cE_1\left(g_1\cE_1-\cA_1(\ri\omega_0)\right)^{-1}\cB u_0\right) \\ \vdots \\ -\operatorname{Re}\left(v_0^\H \cC\left(g_m\cE_m-\cA_m(\ri\omega_0)\right)^{-1}\cE_m\left(g_m\cE_m-\cA_m(\ri\omega_0)\right)^{-1}\cB u_0\right)  \end{bmatrix}.\label{grad formula}
\end{align}
%In the next section we would like to accelerate calculation of $\cH_\infty$ norm by  the using interpolatory reduction techniques derived for particular setting that uses $\cH_\infty$ norm.
\section{Our Approach}
\subsection{Algorithm Description}
\begin{algorithm}[tb]
 \begin{algorithmic}[1]
  % inputs
  \REQUIRE The matrices $M,\,K,\,C_{\rm int},\,B_2,\,E_2,\,H_1$ defining the system \eqref{eq:MDKs} and the associated matrices $\cB$, $\cC$ and the function $\cD(\cdot,\cdot)$ as in \eqref{TF}, some initial parameters $\tg^{(1)},\,\tg^{(2)},\,\ldots,\,\tg^{(\nu)} \in \R^p$
  \ENSURE Optimal gains $g_* = \argminx_{g \in (\R_{\ge 0})^p} \left\| F(g,\cdot) \right\|_{\cH_\infty}$, where $F$ is the closed-loop transfer function in \eqref{TF}.
  \FOR{$i=1,\,2,\,\ldots,\,\nu$}
   	\STATE For given $\tg^{(i)}$ choose initial interpolation points $\ri \tilde{\omega}_{i,1},\,\ldots,\,\ri \tilde{\omega}_{i,k_i} \in \ri\R$.
 % 	\STATE Compute $\|F(\tg^{(i)},.)\|_{\cH_{\infty}}$ using \texttt{linorm\_subsp} and determine the optimal frequency $\omega_i$ \label{frequencyw0}
  \ENDFOR
  \STATE Construct the initial projection matrices
  \begin{align*}
  V_0 &= \begin{bmatrix} \cD(\tg^{(1)},\ri \tilde{\omega}_{1,1})^{-1}\cB,\,\ldots,\,\cD(\tg^{(1)},\ri \tilde{\omega}_{1,k_1})^{-1}\cB,\,\ldots,\,\cD(\tg^{(\nu)},\ri \tilde{\omega}_{\nu,1})^{-1}\cB,\,\ldots,\,\cD(\tg^{(\nu)},\ri \tilde{\omega}_{\nu,k_\nu})^{-1}\cB \end{bmatrix}, \\
  W_0 &= \begin{bmatrix} \cD(\tg^{(1)},\ri \tilde{\omega}_{1,1})^{-\H}\cC^\H,\,\ldots,\,\cD(\tg^{(1)},\ri \tilde{\omega}_{1,k_1})^{-\H}\cC^\H,\,\ldots,\,\cD(\tg^{(\nu)},\ri \tilde{\omega}_{\nu,1})^{-\H}\cC^\H,\,\ldots,\,\cD(\tg^{(\nu)},\ri \tilde{\omega}_{\nu,k_\nu})^{-\H}\cC^\H
  \end{bmatrix}.
  \end{align*}
  %\WHILE{not converged}
  \FOR{$j=1,\,2,\,\ldots,\,r$ (until convergence)}
   \STATE Set $\tilde{\cD}_{j-1}(g,s) := W_{j-1}^\H \cD(g,s) V_{j-1}$, $\tilde{\cB}_{j-1}:=W_{j-1}^\H \cB$, and $\tilde{\cC}_{j-1} := \cC V_{j-1}$ and $\tF_{j-1}(g,s) := \tilde{\cC}_{j-1}\tilde{\cD}_{j-1}(g,s)^{-1}\tilde{\cB}_{j-1}$.
   %\STATE Compute $\hat{g}_j \in (\R_{\ge 0})^p$ and $\hat{\omega}_j \in \R$ such that
   %$\big\|\tilde{F}_j\big(\hat{g}_j,\ri\hat{\omega}_j\big)\big\|_2 = \inf_{g \in (\R_{\ge 0})^p} \sup_{\omega \in \R} \big\| \tilde{F}_j(g,\ri \omega) \big\|_2$.
   \STATE Compute $\hat{g}^{(j)} \in (\R_{\ge 0})^p$ such that
   $\big\|\tilde{F}_{j-1}\big(\hat{g}^{(j)},\cdot)\big\|_{\cL_\infty} = \min_{g \in (\R_{\ge 0})^p} \big\| \tilde{F}_{j-1}(g,\cdot) \big\|_{\cL_\infty}$ (assume that the minimum exists).
   \STATE Compute $\hat{\omega}_{j} \in \R_{\ge 0} \cup \{\infty\}$ such that
   $\big\|F\big(\hat{g}^{(j)},\ri\hat{\omega}_j\big)\big\|_2 = \max_{\omega \in \R_{\ge 0} \cup\{\infty\}} \big\| F\big(\hat{g}^{(j)},\ri \omega\big) \big\|_2$.
  \STATE Update the subspaces and set
  \begin{equation*}
    V_{j} := \begin{bmatrix} V_{j-1} & \cD(\hg^{(j)},\ri\hat{\omega}_j)^{-1}\cB \end{bmatrix}, \quad W_{j} := \begin{bmatrix} W_{j-1} & \cD(\hg^{(j)},\ri\hat{\omega}_j)^{-\H}\cC^\H\end{bmatrix}.
  \end{equation*}
  \STATE Orthonormalize the columns of $V_j$ and $W_j$, respectively.
 % \ENDWHILE
  \ENDFOR
  \STATE Set $g_* := \hg^{(r)}$.
\end{algorithmic}
\caption{Greedy algorithm for semi-active $\cH_\infty$ damping optimization}
\label{alg:algo1}
\end{algorithm}
In this subsection we present our adaptive interpolation approach for damping optimization which has been adapted from \cite{AliBMV18}. The basic algorithm is given in Algorithm~\ref{alg:algo1}. Note that we have formulated the algorithm only for the case $m=\ell$ to simplify its formulation. If $m\neq\ell$, then the updates of the projection matrices must be adapted as in Subsection~\ref{subsec:hinf}. Some words of explanation are in order. To initialize the algorithm, in lines 1--4 initial projection matrices are constructed. This is done by sampling the parameter domain $(\R_{\ge0})^p$ and the imaginary axis appropriately. Since the choice of good sampling points can be a difficult issue, we will give more details about this in Subsection~\ref{IP}.

In lines 5--11, the optimization loop is performed. In line 6, we construct a reduced-order model and its associated transfer function using the projection matrices. Since the reduced-order model is of low dimension, its parameters can be efficiently optimized in line 7. This is a nonconvex and nonsmooth optimization problem. The nonsmoothness arises from the fact that it might well happen, that the function $g \mapsto \left\| F(g,\cdot) \right\|_{\cH_\infty}$ is continuous but \emph{not differentiable} for $g = g_*$. This is the case, when for $F(g_*,\cdot)$ there exist multiple frequencies at which the $\cH_\infty$-norm is attained (i.\,e., Assumption A1 is violated). Luckily, there exist optimization algorithms that are specialized for this kind of problem. Therefore, we make use of %the methods presented in \cite{LewO08,LewO13} which have been implemented in the MATLAB package \texttt{HANSO}\footnote{\url{http://www.cs.nyu.edu/faculty/overton/software/hanso/index.html}} and
the method presented in \cite{CurtisMitchOverton17} which can handle such problems and which has been implemented in the MATLAB package \texttt{GRANSO}\footnote{downloadable from \url{https://gitlab.com/timmitchell/GRANSO}.}. This software only requires a function handle for the evaluation of the mapping $g \mapsto \left\| F(g,\cdot) \right\|_{\cH_\infty}$ and its gradient as discussed in Section~\ref{sec:eval}. Moreover, constraints such as the elementwise nonnegativity of $g$ can be easily incorporated into \texttt{GRANSO}.
Then in line 8 we evaluate the $\cH_\infty$-norm of the original transfer function at the computed optimal parameter value to find a new optimal frequency for sampling and updating the reduced-order model. This is done using the subspace method of \cite{AliBMSV17} for computing the $\cH_\infty$-norm of a large-scale system which is implemented in the MATLAB code \texttt{linorm\_subsp}\footnote{downloadable from \url{http://www.tu-berlin.de/?186267&L=1}}.
%These require differentiability of the objective function for \emph{almost all} points. In particular, it is allowed that the function does not have to be differentiable at the minimizer.
Finally, in line 9 of Algorithm~\ref{alg:algo1} we update the projection matrices according to the optimal parameters we have just obtained. Note that our projection spaces are constructed in a way such that for $i = 1,\,\ldots,\,\nu$, $j=1,\,\ldots,\,r$, and $\ell = 1,\,\ldots,\,p$ we fulfill the \emph{Hermite} interpolation conditions
\begin{align*}
 F(\tg^{(i)}, \ri\tilde{\omega}_{i,1}) &= \tF_r(\tg^{(i)}, \ri\tilde{\omega}_{i,1}),\,\ldots,\,F(\tg^{(i)}, \ri\tilde{\omega}_{i,k_{i}}) = \tF_r(\tg^{(i)}, \ri\tilde{\omega}_{i,k_{i}}), \quad F(\hg^{(j)}, \ri\hat{\omega}_{j}) = \tF_r(\hg^{(j)}, \ri\hat{\omega}_{j}), \\
 \frac{\partial}{\partial \omega}F(\tg^{(i)}, \ri\omega)\Big|_{\omega = \tilde{\omega}_{i,1}} &= \frac{\partial}{\partial \omega}\tF_r(\tg^{(i)}, \ri\omega)\Big|_{\omega = \tilde{\omega}_{i,1}},\,\ldots,\,\frac{\partial}{\partial \omega}F(\tg^{(i)}, \ri\omega)\Big|_{\omega = \tilde{\omega}_{i,k_{i}}} = \frac{\partial}{\partial \omega}\tF_r(\tg^{(i)}, \ri\omega)\Big|_{\omega = \tilde{\omega}_{i,k_{i}}}, \\
 \frac{\partial}{\partial \omega}F(\hg^{(j)}, \ri\omega)\Big|_{\omega = \hat{\omega}_{j}} &= \frac{\partial}{\partial \omega}\tF_r(\hg^{(j)}, \ri\omega)\Big|_{\omega = \hat{\omega}_{j}}, \\
  \frac{\partial}{\partial g_\ell} F(g, \ri\omega_{i,1})\Big|_{g = \tilde{g}^{(i)}} &= \frac{\partial}{\partial g_\ell}\tF_r(g, \ri\omega_{i,1})\Big|_{g = \tilde{g}^{(i)}},\,\ldots,\,\frac{\partial}{\partial g_\ell}F(g, \ri\omega_{i,k_i})\Big|_{g = \tilde{g}^{(i)}} = \frac{\partial}{\partial g_\ell}\tF_r(g, \ri\omega_{i,k_i})\Big|_{g = \tilde{g}^{(i)}}, \\
 \frac{\partial}{\partial g_\ell}F(g, \ri\hat{\omega}_j)\Big|_{g = \hat{g}^{(j)}} &= \frac{\partial}{\partial g_\ell}\tF_r(g, \ri\hat{\omega}_j)\Big|_{g = \hat{g}^{(j)}},
\end{align*}
This also directly leads to Hermite interpolation conditions between the functions $\sigma_{\max}(F(\cdot,\ri\cdot))$ and $\sigma_{\max}(\tilde{F}_r\big(\cdot,\ri\cdot)\big)$. In particular, we have
\begin{equation*}
 \frac{\partial}{\partial \omega}\sigma_{\max}\big(F(\hg^{(j)}, \ri\omega)\big)\Big|_{\omega = \hat{\omega}_{j}} = \frac{\partial}{\partial \omega}\sigma_{\max}\big(\tF_r(\hg^{(j)}, \ri\omega)\big)\Big|_{\omega = \hat{\omega}_{j}} = 0 \quad \forall\,j=1,\,\ldots,\,r.
\end{equation*}
This property is useful to show a superlinear rate of convergence \cite{AliBMV18}, at least in the case of convergence to a differentiable minimizer.
Finally, for avoiding ill-conditioning of the reduced problems, we orthonomalize the columns of $V_j$ and $W_j$, respectively in line 10. 

In this work we consider systems with a larger number of inputs and outputs. In
order to avoid a fast growth of the subspace dimensions in our projection
approach, we make use of \emph{tangential interpolation}. That is, with $v_1
\in \C^{\ell},\,w_1 \in \C^m$ being the right and left singular vectors
associated with the largest singular value of
$F(\hat{g}^{(j)},\ri\hat{\omega}_j)$ we replace
Step~9 in Algorithm~\ref{alg:algo1} by
\begin{equation*}
 V_{j} := \begin{bmatrix} V_{j-1} & \cD(\hg^{(j)},\ri\hat{\omega}_j)^{-1}\cB v_1
\end{bmatrix}, \quad W_{j} := \begin{bmatrix} W_{j-1} &
\cD(\hg^{(j)},\ri\hat{\omega}_j)^{-\H}\cC^\H w_1 \end{bmatrix},
\end{equation*}
with the function $\cD(\cdot,\cdot)$ as in \eqref{TF}. This is similarly done in the projection approach for the computing the
$\cH_\infty$-norm. In this way, the subspace dimension grows by one in each step
which makes these computations feasible even for larger iteration counts. On the
other hand, the convergence analyses of \cite{AliBMSV17,AliBMV18} are not valid
any longer, since we may loose the Hermite interpolation properties between the
maximum singular values functions of the full and the reduced transfer
functions. This problem can be solved by evaluating the full and the reduced
transfer function at the interpolation points, checking whether the Hermite
interpolation properties are satisfied, and expand the projection spaces by
including further singular vectors if necessary. Remarkably, in practice it happens very
rarely that the Hermite interpolation conditions are not satisfied.
%experiments still indicate a very good performance in this case.

Also note that in course of the iteration we use $\ri\hat{\omega}_j$ as a further initial interpolation point for warm-starting the computation oof $\left\|F(g^{j+1},\cdot)\right\|_{\cH_\infty}$ in line~8 of Algorithm~\ref{alg:algo1} for faster convergence. Moreover, the value of $\hat{g}^{(j)}$ is used as an initial value for the optimization of $g \mapsto \big\|\tF_j(g,\cdot) \big\|_{\cL_\infty}$ in line~7 of Algorithm~\ref{alg:algo1}.

\subsection{Choice of the Initial Interpolation Points}
\label{IP}
%To optimize the parameter $g$, we use a gradient descent method. Therefore, the evaluation of the objective function and its gradient needs to be performed numerous times which is computationally very demanding. To reduce the execution time, we will use a greedy subspace method to compute the $\mathcal{H}_\infty$ of large-scale system, see \cite{AliBMSV17}.
Recall that our method for computing the optimal gains $g_*$ is based on Hermite
interpolation between the original and the reduced parameter-dependent
transfer functions. Thus, in case of convergence, we would only know that
\eqref{eq:optimal} is satisfied locally. While convergence to a local minimizer
$g_*$ of $g \mapsto \left\| F(g,\cdot) \right\|_{\cH_\infty}$ would be feasible,
convergence to a local maximizer of $\omega \mapsto \left\| F(g,\ri \omega)
\right\|_{2}$ would be troublesome, since this would underestimate the
closed-loop performance for a given $g$.

For this reason it is important to find good initial points in line 2 of Algorithm~\eqref{alg:algo1} to enhance the chance of converging to a global maximizer of the singular value function.
%However, this method is based on Hermite interpolation and provides only a local maximizer of the largest singular value, not necessarily a global one. To maximize the chances to find the global maximum,
%In our algorithm we very carefully choose the initial interpolation points.
Thus, for $\tg^{(i)},\,\ldots,\,\tg^{(\nu)}$, we estimate the location of the dominant poles of the transfer function which are responsible for large maximum singular values of the transfer functions $F(\tg^{(i)},\cdot)$ on the imaginary axis.

According to \cite{RomM08,BKTT15}, a transfer function $F(s) = H_1(s^2M + sC + K)^{-1}E_2 \in \R(s)^{m \times \ell}$ can be written in pole/residue representation as
\begin{equation*}
 F(s) = \sum_{i=1}^{2n} \frac{R_i}{s-\lambda_i}, %\qquad \text{with} \qquad R_i = (H_1x_i)(y_i^*E_2)\lambda_i
\end{equation*}
where $\lambda_i$ and $R_i$, $i=1,\,\ldots,\,2n$, are the poles and the residues of the transfer function, respectively. Now let
\begin{equation*}
(\lambda_i^2 M + \lambda_i C +  K)x_i = 0 \qquad \text{and} \qquad y_i^\H (\lambda_i^2 M + \lambda_i C +  K) = 0,
\end{equation*}
where $x_i \in \C^n\setminus\{0\}$ and $y_i \in \C^n\setminus\{0\}$ are scaled such that $-y_i^\H K x_i + \lambda_i^2 y_i^\H M x_i = 1$. Then we have that
\begin{equation}\label{eq:res}
 R_i = \lambda_i(H_1x_i)(y_i^\H E_2).
\end{equation}
% With a fixed $g$, the transfer function can be rewritten as
% $$F(s) = \begin{bmatrix} H_1 	& 	0\end{bmatrix}
% 		  \left( s \begin{bmatrix} I_n 	& 	0\\ 	0 	& 	I_n\end{bmatrix}
% 		  - \begin{bmatrix} 0 	& 	I_n\\ 	-M^{-1}K 	& 	-M^{-1}C(g)\end{bmatrix} \right)
% 		  \begin{bmatrix} 0 	\\ 	M^{-1}E_2\end{bmatrix}.$$
% According to \cite{BKTT15}, since the transfer function is meromorphic, it holds
% \bea
% \label{sum}
% F(s) = \sum_{i=1}^{2n} \frac{R_i}{s-\lambda_i} \qquad \text{with} \qquad R_i = (H_1x_i)(y_i^*E_2)\lambda_i,
% \eea
% where $\lambda_i \in \C$ are the eigenvalues of the corresponding quadratic eigenvalues problem:
% $$(\lambda_i^2 M + \lambda_i C(g) +  K)x_i = 0 \qquad \text{and} \qquad y_i^*(\lambda_i^2 M + \lambda_i C +  K) = 0,$$
% with $x_i, y_i \in \C^{2n} \backslash \{ 0 \}$, respectively right and left eigenvectors.

%For a proper approximation of the transfer function $F(s)$, we want to compute only the significant terms in the sum. Definition of significant terms depends on considered applications and we will focus the the one for which it is shown that plays important role in damping optimization setting, for more details, see, e.g. \cite{BKTT15,TomljGugerBeatt17}.

\begin{mydef}\label{def: dominant pole}
For the transfer function $F(s) = H_1(s^2M + sC + K)^{-1}E_2 \in \R(s)^{m \times \ell}$, a pole $\lambda_i$ is called \emph{dominant} if
\begin{equation*}
  \frac{\left\|R_i\right\|_2}{|\Real{\lambda_i}|} > \frac{\left\|R_j\right\|_2}{|\Real{\lambda_j}|} \quad \forall\, j \neq i
\end{equation*}
The $k$ most dominant poles are the $k$ poles with the largest values of $\frac{\left\|R_i\right\|_2}{|\Real{\lambda_i}|}$.
\end{mydef}
The dominant poles have a close relationship to the local maxima of the function
$\omega \mapsto \left\| F(\ri\omega) \right\|_2$. Namely, if $\lambda_i$ is a
dominant pole of $F(s)$, then it is likely that a large local maximum can be
found near $\operatorname{Im}(\lambda_i)$. Conversely, if $\omega \mapsto
\left\| F(\ri\omega) \right\|_2$ has a large local maximum at $\omega_0$, then
a dominant pole will likely be close to $\ri \omega_0$. A graphical
interpretation of this relationship is available in \cite[p.
17]{Kur10}. A similar notion of pole dominance has already been used in the
context of damping optimization, see \cite{BKTT15}.

In principal, one could compute a number of dominant poles of $F(\tg^{(i)},s) \in \R(s)^{m \times \ell}$ for all $\tg^{(i)}$, $i=1,\,\ldots,\,\nu$ and use their imaginary parts as initial frequencies $\tilde{\omega}_{i,j}$ in line 2 of Algorithm~\ref{alg:algo1}. This can be done using the subspace accelerated MIMO dominant pole algorithm (SAMDP) \cite{RomM08}.
%
%for second-order systems (SAQMDP) which is outlined in \cite[Algorithm 1]{BKTT15}.

%The problem is to estimate $R_i$ and $\lambda_i$ for the damped problem in order to determine the dominant poles. Then we will use them as starting points for the evaluation of the $\mathcal{H}_\infty$ norm of the transfer function.

%One approach for efficient calculation of dominant poles within considered mechanical system was described in \cite{BKTT15,TomljGugerBeatt17}, but in order to have initial shifts that can be calculated even more efficiently we propose the following approach.

On the other hand, if the state-space dimension $n$ is of moderate size, then we can also estimate $\left\|R_i\right\|_2$ and $\lambda_i$ by transforming the system to modal coordinates before. This will be outlined next.
We estimate $\left\|R_i\right\|_2$ using the transfer function of the undamped problem which is
\begin{equation*}
F_0(s) = H_1(s^2M+K)^{-1}E_2.
\end{equation*}
Since $M$ and $K$ are positive definite matrices, there exists a matrix $\Phi = \begin{bmatrix} \phi_1 & \ldots & \phi_n \end{bmatrix}\in \R^{n \times n}$ such that
\begin{equation*}
\Phi^\T(s^2M+K)\Phi = s^2I_n+\Omega^2, \quad \text{where} \quad \Omega = \diag(\omega_1,\ldots,\omega_n), \quad \omega_1 > \omega_2 > \ldots > \omega_n > 0.
\end{equation*}
Therefore, we can write
\begin{equation*}
F_0(s) = (H_1\Phi)(s^2I_n+\Omega^2)^{-1}(\Phi^\T E_2).
\end{equation*}
Then by \eqref{eq:res}, for $i=1,\,\ldots,\,n$, the residues of $F_0$ are given by
\begin{equation*}
 R_{0,i}^\pm = \pm \ri\gamma_i\omega_i(H_1\Phi e_i)(e_i^\T \Phi^\T E_2) = \pm \ri\gamma_i\omega_i(H_1\phi_i)(\phi_i^\T E_2),
\end{equation*}
where $e_i \in \Rn$ is the $i$th unit vector in $\Rn$ and $\gamma_i$ is a scaling parameter such that
\begin{equation*}
 \gamma_i \big(-e_i^\T \Omega^2 e_i - \omega_i^2 e_i^\T e_i\big) = 1.
\end{equation*}
This gives $\gamma_i = -\frac{1}{2\omega_i^2}$ and hence, we get
\begin{equation*}
 \left\|R_{0,i}^\pm\right\|_2 = \frac{1}{2}\omega_i^{-1} \left\|H_1\phi_i\right\|_2 \left\|E_2^\T \phi_i\right\|_2.
\end{equation*}

Next we estimate the pole locations of the transfer function $F(s)$. By a linearization of the quadratic eigenvalue problem
\begin{equation*}
 \big(\lambda^2 I_n + \lambda \Phi^\T C(g) \Phi + \Omega^2\big) x = 0 \quad \text{with} \quad x \in \C^n \setminus \{0\},
\end{equation*}
the poles of $F(s)$ are eigenvalues of the matrix
\begin{equation}
 \cA := \begin{bmatrix} 0 & \Omega \\ -\Omega & -\Phi^\T C(g) \Phi \end{bmatrix}.
\end{equation}
Here we consider the matrix $\cA$ as a perturbation of the matrix
\begin{equation}
 \cA_0 := \begin{bmatrix} 0 & \Omega \\ -\Omega & 0 \end{bmatrix}.
\end{equation}
%We consider the matrices
%$$A = \begin{bmatrix} 0 & \Omega \\ -\Omega & 0\end{bmatrix} \qquad \text{and }\Delta A = \begin{bmatrix} 0 & 0 \\ 0 & \Phi^T(C_{\rm int}+B_2GB_2^T)\Phi\end{bmatrix}$$
%where the matrix $\Delta A$ represents the perturbation of the damping to the undamped matrix $A$.
Note that the eigenvalues of $A_0$ are $\lambda_{0,i}^\pm = \pm \ri \omega_i$ with associated right and left normalized eigenvectors (both are the same)
\begin{equation*}
 v_{0,i}^\pm = \frac{1}{\sqrt{2}} \begin{pmatrix} e_i \\ \pm \ri e_i \end{pmatrix}.
\end{equation*}
By a standard result from first order perturbation theory \cite{SteS90}, it holds that the eigenvalues $\lambda_i^{\pm}$, $i = 1,\,\ldots,\,n$ of $\cA$ are given by
\begin{align*}
   \lambda_i^{\pm} &= \pm \ri \omega_i + (v_{0,i}^\pm)^\H \begin{bmatrix} 0 & 0 \\ 0 & -\Phi^\T C(g) \Phi \end{bmatrix} v_{0,i}^\pm + \cO\left(\left\|C(g)\right\|_2^2\right) \\
   &= \pm \ri \omega_i + \frac{1}{2} \phi_i^\T C(g) \phi_i + \cO\left(\left\|C(g)\right\|_2^2\right),
\end{align*}
which gives the approximation
\begin{equation*}
 \left|\Real{\lambda_i^\pm}\right| \approx \frac{1}{2} \phi_i^\T C(g) \phi_i.
\end{equation*}
This approximation is usually good, if the norm of the perturbation given by $\left\|C(g)\right\|_2$ is small.
Finally we sort $\omega_i$ according to the values of
\begin{equation}\label{aprox real part}
 \operatorname{dom}(\omega_i) :=
\frac{\left\|R_{0,i}^+\right\|_2}{\left|\Real{\lambda_i^+}\right|}
\end{equation}
in decreasing order, which can be done very cheaply once the system is known in modal coordinates.

% where $\lambda_i$ is the perturbed eigenvalue.
% In our case, $\lambda_i$ is purely imaginary, so we have
% \begin{eqnarray}
% \Real{\lambda_i} & = & \text{Re} (\ri \omega_i) + \text{Re} \left( \begin{bmatrix}  -\ri \phi_i & \phi_i \end{bmatrix}^T \begin{bmatrix} 0 & 0 \\ 0& \Phi^T(C_{\rm int}+B_2GB_2^T)\Phi \end{bmatrix} \begin{bmatrix} \ri \phi_i \\ \phi_i  \end{bmatrix} \right) \\
% 		         & = &  \phi_i^T \cdot \Phi^T(C_{\rm int}+B_2GB_2^T)\Phi \cdot \phi_i . \label{aprox real part}
% \end{eqnarray}

%The above approach can be used for determination of initial points. Derived approximations of eigenvalues that can be used for estimation of dominant poles will be accurate for small initial gains. However, in numerical experiments we have observed that these initial points can provide efficient calculation of $\cH_\infty$ norm.

%Alternatively, dominant poles defined in \eqref{def: dominant pole}  can be calculated by subspace accelerated quadratic MIMO dominant pole algorithm (SAQMDP). This algorithm can be implemented such that it uses second order structure efficiently, for more details see, e.g. \cite{BKTT15, RomM06a}.

%Once when we have calculated $k$ dominant poles we can use imaginary parts as the initial points for the calculation of  $\cH_\infty$ norm. In particular, in this paper we will compare usage of initial points based on above mentioned approach that uses first order approximation with initial points that are calculated from dominant poles calculated by Algorithm 1 SAQMDP from  \cite{BKTT15}.

%with a tolerance of $10^{-16}$ for (approximate) stationarity

\section{Numerical Experiments}
In this section we consider numerical experiments. In order to illustrate the efficiency of our proposed approach, we compare the optimal gains calculated for the full-order system with the optimal gains obtained by Algorithm \ref{alg:algo1}. In each of our experiments we use the open-source package \texttt{GRANSO} for the parameter optimization. Alternatively, one could use \texttt{HANSO}\footnote{available at \url{http://www.cs.nyu.edu/faculty/overton/software/hanso/index.html}}, see also \cite{BurLO05,LewO08}. However, since we wish the optimal parameters to be in a feasible domain, we have a constrained optimization problem which makes the use \texttt{GRANSO} more appropriate in our situation.
For the full-order problem, the $\cH_\infty$-norms are calculated by the Boyd-Balakrishnan algorithm for descriptor systems \cite{BenSV12} which is implemented in the FORTRAN routine \texttt{AB13HD.f} and is called by the mex file \texttt{linorm\_h.f} in \texttt{MATLAB}.
On the other hand, in Algorithm~\ref{alg:algo1}, we compute the $\cH_\infty$-norms for a sequence of large-scale non-parametric problems by employing \texttt{linorm\_subsp} and then optimize the gains on a sequence of reduced parametric systems.
%In all tests, the gradient of the objective function is calculated with the formula \eqref{grad formula}.

%
In our experiments, the computations have been carried out on a machine with
four Intel{\textregistered} Core{\texttrademark}
i5-4590 CPUs @ 3.30\,GHz and 16\,GB RAM. The results
reported in this work are calculated by {\tt MATLAB} version 9.6.0.1072779
(R2019a) on a 64-bit Linux operating system.

%Intel\textsuperscript{\textregistered} Xeon\textsuperscript{\textregistered}
%CPU E5-1650 v2 @ 3.50\,GHz, 12 CPUs, and 16\,GB RAM. The results reported in
%this work are calculated by {\tt MATLAB} version 8.3.0.532 (R2014a) on a
%64-bit Linux operating system.

\begin{example}
We consider an $n$-mass oscillator which
describes the mechanical system of $n$ masses and $n+1$ springs shown in Figure \ref{nmass_osc_semi_active}.
The mathematical model is given by
\eqref{MDK1s}--\eqref{MDK3s}, while the mass and stiffness matrices are
\begin{align*}
M & = \diag(m_1,m_2,\ldots,m_n)  ,\\
K &=
\begin{bmatrix}
  k_1+k_2 & -k_2 &  &  &  \\
  -k_2    & k_2+k_3& -k_3 &  &  \\
          & \ddots  & \ddots & \ddots &  \\
          &      & -k_{n-1} & k_{n-1}+k_{n} & -k_n   \\
          &      &  & -k_n & k_n+k_{n+1} \\
\end{bmatrix}.
\end{align*}
The masses and stiffnesses have the following configuration:
\begin{equation*}
 n = 700; \quad k_i = 10, \quad i=1,\,\ldots,\,700; \quad
   m_i  =  \begin{cases}  200.3 - 0.6 i, \quad & i=1,\,\ldots,\,300, \\
      0.4 i-100.2, \quad & i=301,\,\ldots,\,700.
                 \end{cases}
\end{equation*}

\begin{figure}[tb]
\centering
\includegraphics[width=10cm]{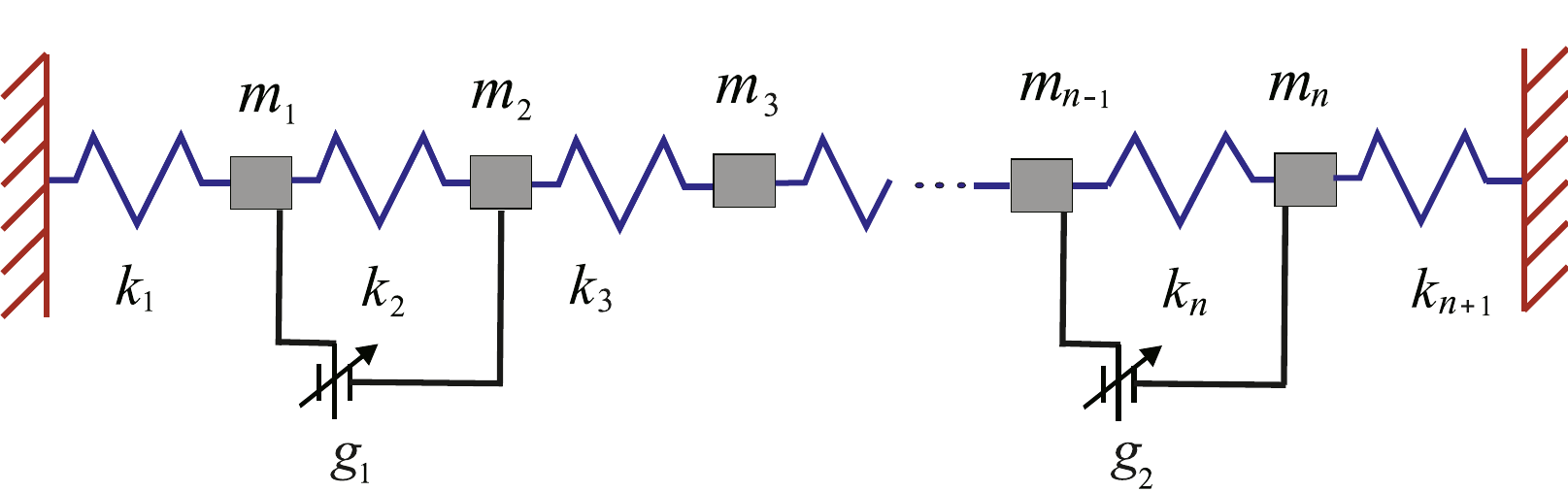}
\caption{An $n$-mass oscillator with two dampers}
\label{nmass_osc_semi_active}
\end{figure}

We are interested in the states that correspond to the masses with indices ranging from 290 to 309, that is, in this example we consider the performance output
 \[
 z(t) = \begin{pmatrix} q_{290}(t) \\ q_{291}(t) \\ \vdots \\ q_{309}(t) \end{pmatrix}.
 \]
Hence, we choose $H_1\in\mathbb{R}^{20\times n}$ where
\begin{equation*}
 H_1(1:20,290:309) = I_{20\times 20}
\end{equation*}
and all other entries are equal to zero.

We define $E_2\in \mathbb{R}^{n\times 10}$ that corresponds to the excitation of the masses that are closer to the grounded masses. More precisely, we choose the matrix $E_2$ to be zero everywhere except for
\begin{equation*}
  E_2(1:5,1:5) = \diag(5,4,3,2,1) ,\quad E_2(696:700,6:10) = \diag(5,4,3,2,1) .
\end{equation*}

The damping matrix $D$ is equal to $C_{\rm int}+B_2GB_2^\T$, where the internal damping $C_{\rm int}$ is given by \eqref{Cint} and where we choose the critical damping model
\begin{equation*}
C_{\rm crit} = 2M^{1/2}\left(M^{-1/2}KM^{-1/2}\right)^{1/2}M^{1/2} > 0,
\end{equation*}
which is widely used in the literature, see, e.\,g., \cite{BRAB98, BennerTomljTruh11, BKTT15}.
Note that this choice makes the unforced system \eqref{MDK1s} asymptotically stable.%, but on the other hand, the matrices than define closed-loop system will in general be dense, since $C_{\rm crit}$ is dense. Thus, the use of sparse arithmetics will not be of great advantage here.

We consider two dampers with different gains, that is, we choose the matrix  $G=\diag(g_1,g_2)\in \mathbb{R}^{2\times 2}$. The geometry of the external damping depends on the dampers' positions $(j,k)$ which are encoded in the matrix $B_2$ by setting
\begin{equation*}
 B_2  = \begin{bmatrix} e_{j}-e_{j+1} &  e_{k}-e_{k+1}\end{bmatrix},
\end{equation*}
 where $e_j$ and $e_k$ are the $j$th and the $k$th canonical vector in $\R^n$.
 %Here, we have introduced an indices $j$ and $k$ since we will consider different damping positions because our main aim is to determine optimal gains for the set of different damping positions.
\end{example}

In general, the problem is to optimize the positions of the dampers, but this requires an optimization of the gains for a large number of different damper positions. Thus, we illustrate the efficiency of our method for the optimization of the damper positions $(j,k)$ for $j \in\{ 40,\, 140,\, 240,\, 340,\, 440,\, 540 \}$ and $k\in \{60,\, 160,\, 260,\, 360,\, 460,\, 560\}$, so in total we obtain 36 different configurations for which we have to perform a gain optimization.

Moreover, the parameter $\alpha_c$ determines the influence of an internal damping as shown in \eqref{Cint}. The internal damping has a strong influence on the system, so we show the results of our approach for two different settings:
\begin{enumerate}[\bf \mbox{Problem} a)]
 \item $\alpha_c = 10^{-5}$: In this case we have a very small internal damping which means that almost all eigenvalues of the matrix polynomial $s^2M + sC(0) + K$ (that corresponds to the problem without dampers) are very close to the imaginary axis. The optimization of this problem is particularly difficult, since these eigenvalues often introduce a very large number of thin peaks (local maxima) in the function $\omega \mapsto \sigma_{\max}(F(g,\ri\omega))$.
 \item $\alpha_c = 10^{-2}$: In this case the internal damping is moderate which improves asymptotic stability. Moreover, the amount of large local maxima in $\omega \mapsto \sigma_{\max}(F(g,\ri\omega))$ is moderate which makes it easier to compute the $\cH_\infty$-norm in this case. This choice of $\alpha_c$ is more realistic from the practitioner's point of view.
\end{enumerate}
For both cases we use the following general computational setup and parameters in Algorithm~\ref{alg:algo1}:
\begin{itemize}
 \item The tolerance of \texttt{GRANSO} for (approximate) stationarity is set to $10^{-12}$.
 \item We use \texttt{linorm\_subsp} with default parameters and with tangential interpolation. For each call of \texttt{linorm\_subsp}, the 30 most dominant poles computed according to Subsection~\ref{IP} are chosen as initial interpolation points. In case we use SAMDP, we first compute 180 dominant poles and pick the 30 most dominant ones out of these.
 \item Initially, we choose 4 initial parameters to set up the initial reduced-order model as follows.
\begin{itemize}
\item For \textbf{Problem a)} the starting initial parameters are $\tilde{g}^{(1)} = \left(\begin{smallmatrix} 10 \\ 10\end{smallmatrix}\right)$, $\tilde{g}^{(2)} = \left(\begin{smallmatrix} 10 \\ 100\end{smallmatrix}\right)$, $\tilde{g}^{(3)}=\left(\begin{smallmatrix} 100 \\ 10\end{smallmatrix}\right)$, and $\tilde{g}^{(4)} = \left(\begin{smallmatrix} 100 \\ 100\end{smallmatrix}\right)$.
\item For \textbf{Problem b)} the starting initial parameters are chosen as $\tilde{g}^{(1)} = \left(\begin{smallmatrix} 100 \\ 100\end{smallmatrix}\right)$, $\tilde{g}^{(2)} = \left(\begin{smallmatrix} 100 \\ 1000\end{smallmatrix}\right)$, $\tilde{g}^{(3)}=\left(\begin{smallmatrix} 1000 \\ 100\end{smallmatrix}\right)$, and $\tilde{g}^{(4)} = \left(\begin{smallmatrix} 1000 \\ 1000\end{smallmatrix}\right)$.
\end{itemize}
 \item The relative termination tolerances for the gains and the $\cL_\infty$-norm are both set to $10^{-6}$. In other words, we terminate, if
 \begin{align*}
  \left\| \hat{g}^{(j)} - \hat{g}^{(j-1)} \right\|_2 &< \frac{1}{2}\cdot 10^{-6} \cdot\left\| \hat{g}^{(j)} + \hat{g}^{(j-1)} \right\|_2, \quad \text{or} \\
  \left| \big\|\tilde{F}_j(\hat{g}^{(j)},\cdot)\big\|_{\cH_\infty} - \big\|\tilde{F}_{j-1}(\hat{g}^{(j-1)},\cdot)\big\|_{\cH_\infty} \right| &< \frac{1}{2}\cdot 10^{-6} \cdot \left| \big\|\tilde{F}_j(\hat{g}^{(j)},\cdot)\big\|_{\cH_\infty} + \big\|\tilde{F}_{j-1}(\hat{g}^{(j-1)},\cdot)\big\|_{\cH_\infty} \right|,
 \end{align*}
 or the maximum number of iterations (which we have set to 30) has been reached.
\end{itemize}

We further test our method with different computational modes whose results we present further below:
\begin{enumerate}[\bf \mbox{Mode} i)]
 \item For each $\tilde{g}^{(1)},\,\ldots,\,\tilde{g}^{(4)}$, we perform tangential interpolation of $F(\tilde{g}^{(i)},\cdot)$ at 30 equidistantly distributed samples $\ri\tilde{\omega}_{i,1},\,\ldots,\,\ri\tilde{\omega}_{i,30}$ in the frequency range $\tilde{\omega}_{i,j} \in [0,\omega_{\max}]$ for $j=1,\,\ldots,\,30$  to setup the initial reduced-order model. Here $\omega_{\max}$ is the maximum modulus among all eigenvalues of the matrix polynomial $s^2M + K$. To determine initial interpolation points for \texttt{linorm\_subsp} we choose the heuristic approach from Section~\ref{IP}.
 \item We construct the initial reduced-order model as in \textbf{Mode i)}. To determine initial interpolation points for \texttt{linorm\_subsp} we use the SAMPD algorithm.
 \item We compute the $\cH_\infty$-norms of $F(\tilde{g}^{(i)},\cdot)$ for $i = 1,\,2,\,3,\,4$. If the norms are attained at $\ri\tilde{\omega}_1,\,\ldots,\,\ri\tilde{\omega}_4$, we do tangential interpolation only at $\big(\tilde{g}^{(1)},\ri\tilde{\omega}_1\big),\,\ldots,\,\big(\tilde{g}^{(4)},\ri\tilde{\omega}_4\big)$ to construct the initial reduced-order model. To determine initial interpolation points for \texttt{linorm\_subsp} we choose the heuristic approach from Section~\ref{IP}.
\item We construct the initial reduced-order model as in \textbf{Mode iii)}. To determine initial interpolation points for \texttt{linorm\_subsp} we use the SAMDP algorithm.
\end{enumerate}

%
%  Thus, we will present our comparison for two cases, that is, for the case where the internal damping is very small and to the case where internal damping is moderate (in percentages compared to the critical damping).
% In both cases within the application of algorithm \texttt{linorm\_subsp} we have used tolerance $1\cdot 10^{-6}$ and 80 initial interpolation points in Step 2 of Algorithm \ref{alg:algo1} (i.e., $k_i=80$, $i=1,\ldots,\nu$). Also, in both cases we compare different strategies for determination of initial interpolation points.

%Thus, we will present results obtained  by  new approach that is based on formula \eqref{aprox real part} and approach that is based on Algorithm SAQMDP from \cite[Algorithm 1]{BKTT15}. The initial value for the frequency  in SAQMDP is taken as the eigenvalue closest to the imaginary axis corresponding to the zero initial gain.

%\textbf{Case I)} First, we will consider the case where internal damping given by \eqref{Cint} is smaller, in particular, we will have that $\alpha_c=0.00001$.

% In order to present relative errors more clearly, the different configurations are sorted w.r.t. the magnitude of the relative error in the optimal gain.

In the following we compare our results to the naive optimization approach that consists of optimizing the full-order problem with \texttt{GRANSO}
% \todo{Which termination tolerance? The tolerance for approximate stationarity \texttt{opt\_tol}?, ZT: yes, it was 10^{-4}}
and the initial point set to be the first initial parameter, that is the initial point for \textbf{Problem a)} is $\left(\begin{smallmatrix} 10 \\ 10\end{smallmatrix}\right)$, while for \textbf{Problem b)} it is set to $\left(\begin{smallmatrix} 100 \\ 100\end{smallmatrix}\right)$.

%Additionally, in the Algorithm \ref{alg:algo1} we have used three initial parameters that is, parameter $\nu=3$ with initial parameters $\tg_1=[10,\,10]$, $\tg_2=[100,\,10]$ and $\tg_3=[10,\,100]$.

%In the following we illustrate the behavior of the proposed algorithm for viscosity optimization with the setups discussed above. In particular, we consider the relative errors between the optimal gains and the optimal $\mathcal{H}_\infty$-norms

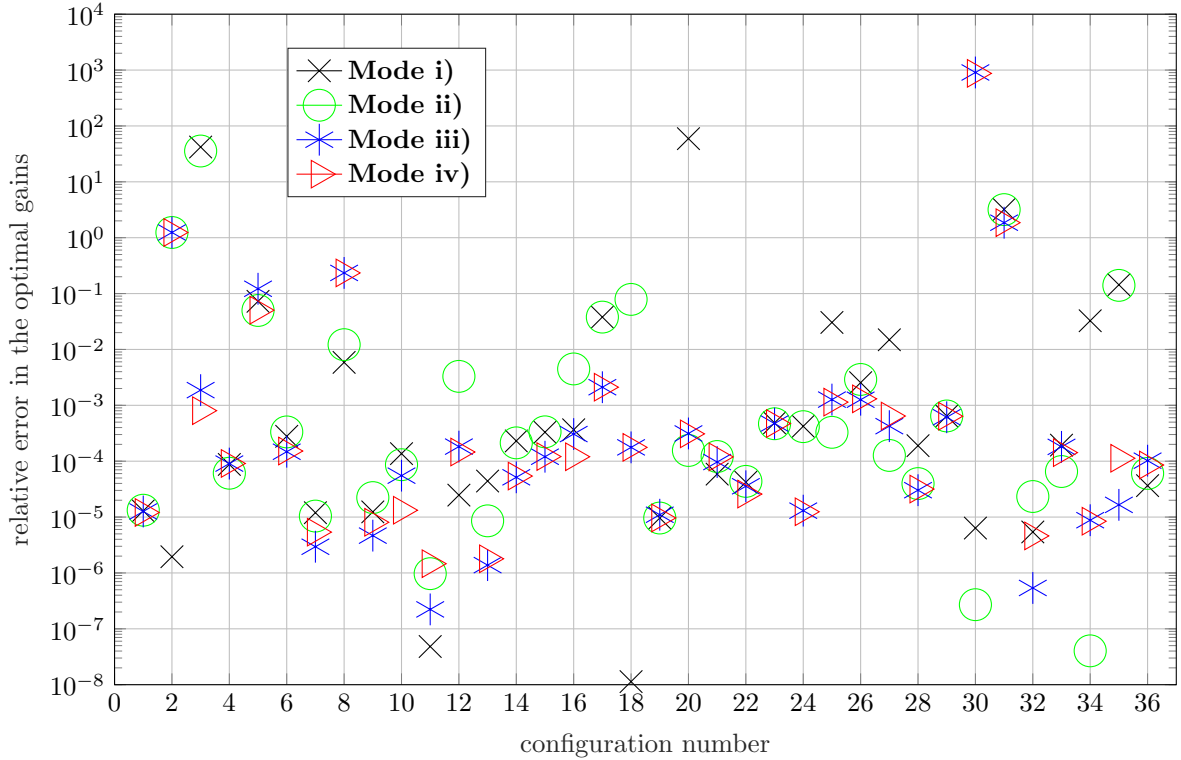
\begin{figure}[tb]\begin{center}
 \input{relErrviscAproblemUpdated.tikz}
 \caption{Comparison of optimal gains for \textbf{Problem a)}} \label{opt gains}\end{center}
\end{figure}

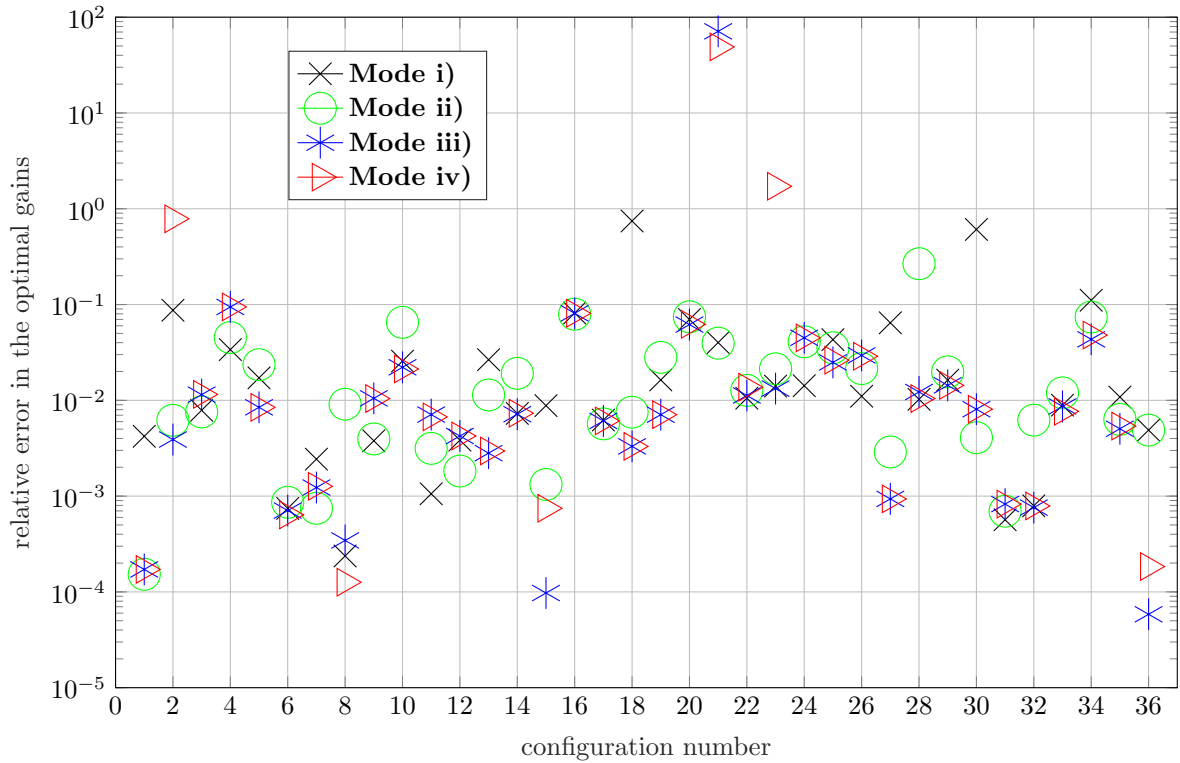
\begin{figure}[tb]\begin{center}
 \input{relErrviscBproblemUpdated.tikz}
\caption{Comparison of optimal gains for \textbf{Problem b)}} \label{opt gainsII}\end{center}
\end{figure}

Figures~\ref{opt gains} and~\ref{opt gainsII} show the relative errors in the computed optimal gains for \textbf{Problem a)} and \textbf{Problem~b)}, respectively, where as reference values we choose the optimal gains computed for the full-order model. More precisely, for the $i$th configuration, the Figures display the values of ${\|g_*^{(i)} -\tilde g_*^{(i)} \|}_2/{\|g_*^{(i)} \|}_2$, where $g_*^{(i)}$ and $\tilde g_*^{(i)}$ denote the gains obtained by optimizing the full-order model directly and by applying Algorithm~\ref{alg:algo1}, respectively.

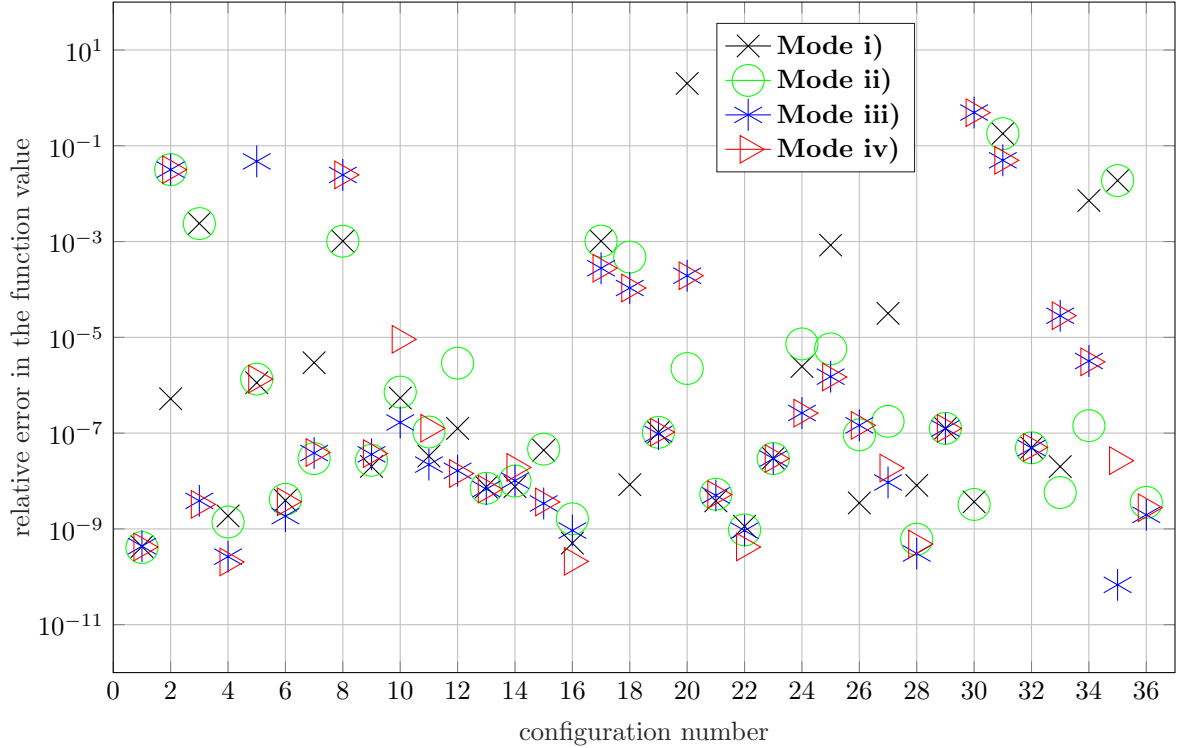
\begin{figure}[tb]\begin{center}
 \input{relErrHinfAproblemUpdated.tikz}
\caption{Comparison of the optimal $\mathcal{H}_\infty$-norm for \textbf{Problem a)}} \label{opt_Hinf}\end{center}
\end{figure}

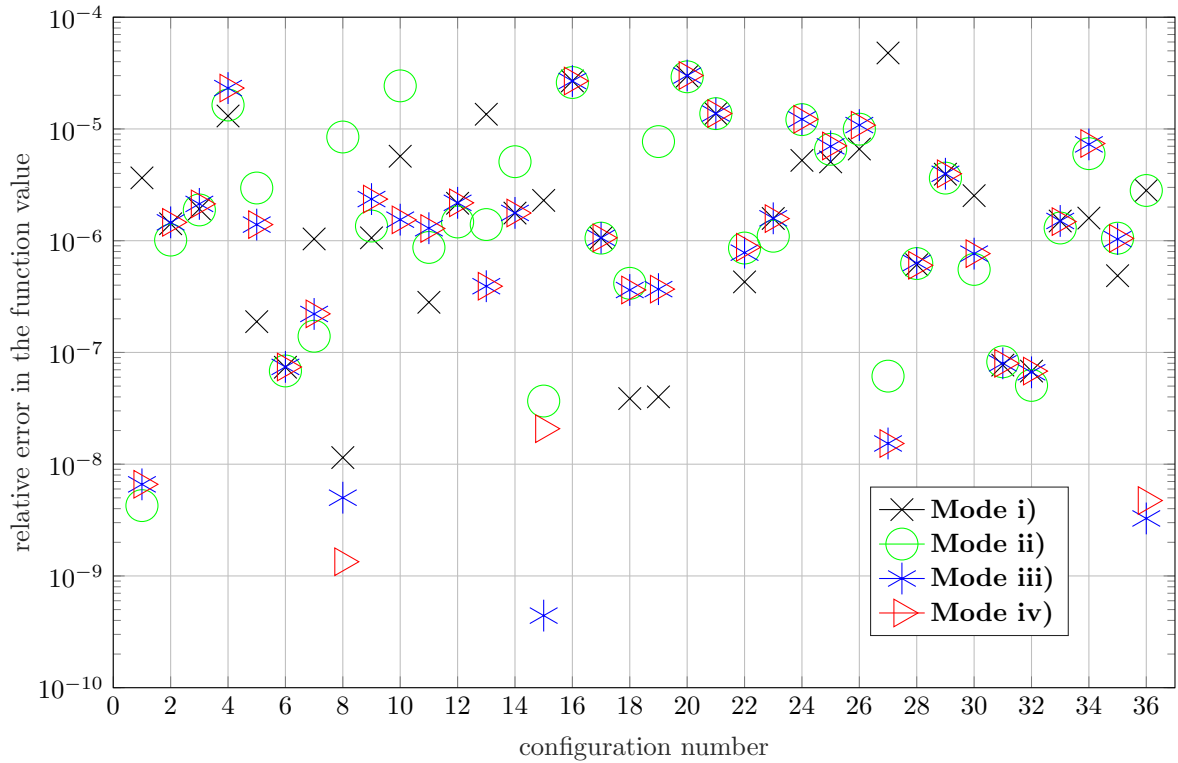
\begin{figure}[tb]\begin{center}
 \input{relErrHinfBproblemUpdated.tikz}
\caption{Comparison of the optimal $\mathcal{H}_\infty$-norm for \textbf{Problem b)}} \label{opt_HinfII}\end{center}
\end{figure}

Figures~\ref{opt_Hinf} and~\ref{opt_HinfII} on the other hand show the relative errors in the computed $\cH_\infty$-norms. In other words, the Figures present the values of ${|{\| F(g_*^{(i)}) \|}_{\mathcal{H}_\infty} - {\| F( \tilde g_*^{(i)}) \|}_{\mathcal{H}_\infty} |}/{{\| F(g_*^{(i)}) \|}_{\mathcal{H}_\infty} }$ for the $i$th configuration. Since the objective function may be flat near the optimizer, i.\,e., there may be a large area of almost optimal gains, the  differences in the optimal values of the objective function are more informative when assessing the quality of the results.

%On this and on the next three succeeding figures \textbf{Mode i)}, \textbf{Mode ii)}, \textbf{Mode iii)} and \textbf{Mode iv)} are represented by black crosses, green circles, blue asterisks and red triangles, respectively.

\begin{figure}[tb]\begin{center}
 \input{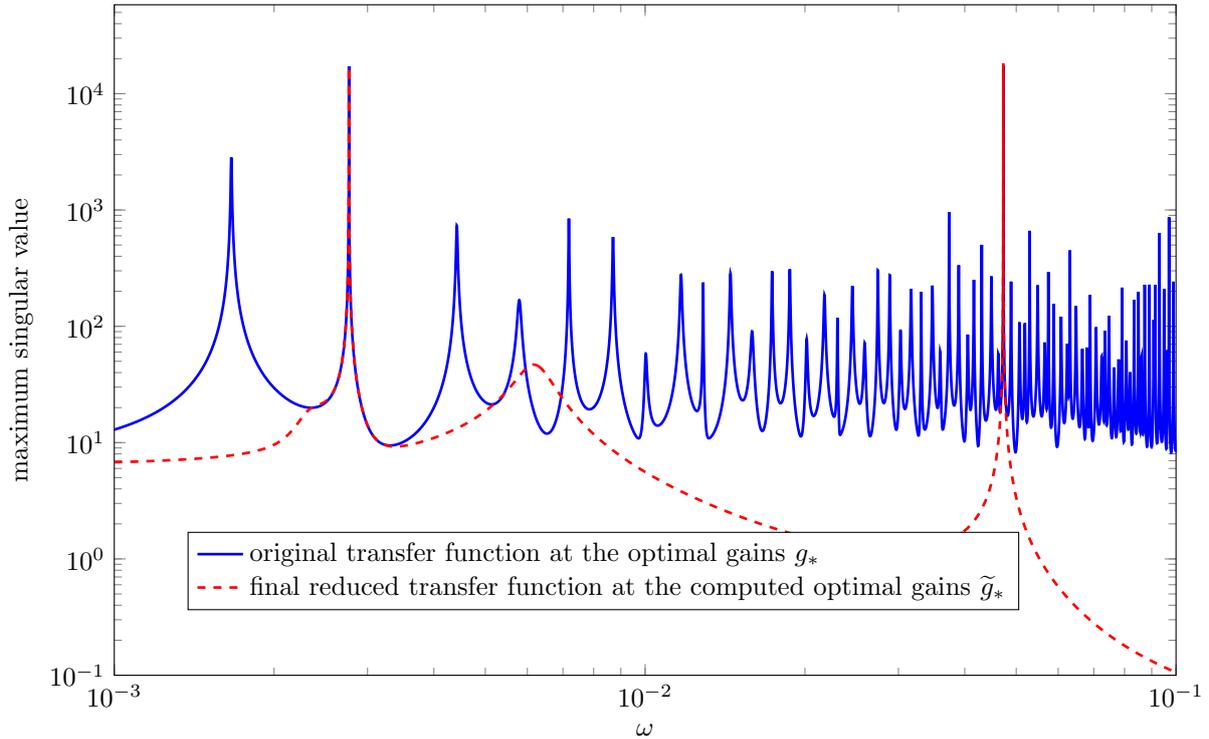}
\caption{Maximum singular value plots for the original transfer function at the optimal gains $g_*$ and the final reduced transfer function at the computed optimal gains $\widetilde{g}_*$ for the 8th configuration of \textbf{Problem~a)} using \textbf{Mode iv)}.} \label{fig:compare}\end{center}
\end{figure}

As expected, our method has some difficulties with \textbf{Problem a)}, since for a very small internal damping, the mappings $\omega \mapsto \left\| F(g,\ri\omega) \right\|_2$ have a huge number of local maxima and then it may e.\,g. happen that the global maximum is missed by our method. Further, these local maxima often lead to extremely steep peaks in the maximum singular value plots which can cause additional numerical difficulties since then the maximum singular values of the transfer function are very sensitive with respect to small changes in $\omega$. These problems can be seen in Figure~\ref{opt_Hinf}, where some of the errors are quite large. This behavior is a clear limitation of Algorithm~\ref{alg:algo1}. The difficulty of this problem is further illustrated in Figure~\ref{fig:compare} in which we plot $\left\|F(g_*,\ri\omega)\right\|_2$ as well as $\big\|\widetilde{F}_r\big(\tilde{g}_*,\ri\omega\big)\big\|_2$ over $\omega$ (where $\widetilde{F}_r$ denotes the final reduced transfer function and $\widetilde{g}_*$ are the optimal gains computed by Algorithm~\ref{alg:algo1}.) for one of the configurations of \textbf{Problem a)}. Remarkably, even if the original transfer function induces extremely many peaks in its maximum singular value plot, this is not the case for the reduced one, but the two global maximizers still coincide approximately.

On the other hand, as shown in Figure~\ref{opt_HinfII}, the results are more convincing for a moderate internal damping which is also more realistic from the application's point of view.
%From presented figures we can conclude  that for the same system configuration, but with different influence of internal damping, our approach is more accurate in the problem where internal damping is moderate. This is natural since that internal damping has strong influence on system stability and robustness of optimal parameters.

We also compare the computational time needed for the optimization process. For
each configuration we measure the runtimes needed to optimize the full-order
model with the naive approach and the time needed to run
Algorithm~\ref{alg:algo1} with the different modes described above. The
corresponding average speedup factors are listed in Table~\ref{time ration}.
\begin {table}[tb]
\caption{Average time ratio} \label{time ration}
\begin{center}
\begin{tabular}{c|cc}
  %\hline
  % after \\: \hline or \cline{col1-col2} \cline{col3-col4} ...
    & \textbf{Problem a)}  & \textbf{Problem b)} \\  \hline
  \textbf{Mode i)}   &  73 &  103 \\
  \textbf{Mode ii)} &   35 &  53\\
  \textbf{Mode iii)} &  295&  434 \\
  \textbf{Mode iv)} &   48 &  67 \\ \hline
\end{tabular}
\end{center}
\end{table}
From this table we can conclude that Algorithm~\ref{alg:algo1} accelerates the
optimization process considerably. For our examples, the initial interpolation
points can be calculated faster by using formula \eqref{aprox real part}
compared to using the SAMDP algorithm, i.\,e., see the entries
corresponding to \textbf{Mode i)} and
\textbf{Mode iii)} in Table \ref{time ration}. However, formula \eqref{aprox real part} can no longer be
efficiently applied, if the system dimension becomes very large. On the other
hand, the SAMDP algorithm could still potentially be used in this case, since
only a small fraction of the eigenvalues may be computed. For this, it is
necessary to be able to solve the involved linear systems efficiently. For
instance, this is the case if instead of critical damping one uses Rayleigh
damping, since then $C_{\rm int} = \alpha M + \beta K$ for some small constants
$\alpha,\,\beta > 0$ would still be a sparse matrix. We do not consider this
case here, since for very large systems, the computations for the naive
optimization approach become prohibitively expensive. From the figures we can
further see that it is more efficient to evaluate the $\cH_\infty$-norms at the
initial gains in order to set up the initial reduced order model (\textbf{Mode
iii)} and \textbf{Mode iv)}) compared the sampling procedure used in
\textbf{Mode i)} and \textbf{Mode ii)}.

We have also evaluated our approach with full instead of tangential
interpolation. In our experience, this does not have a significant influence on
the number of iterations but evaluating the $\cH_\infty$-norm of the reduced
transfer functions becomes very expensive, since the projection space
dimensions grow very fast. Therefore, in this case the runtime is in the same
order of magnitude as for the naive method.

\section{Conclusion}
Altogether we can conclude that with our new approach we have been able to perform
the semi-active $\mathcal{H}_\infty$ damping optimization for a problem with
moderate internal damping with satisfactory relative accuracy, while the
optimization process was considerably accelerated. On the other hand, our method
has a few problems with some of the configurations with very small internal
damping. Such problems must be treated more carefully and it is necessary to
have a
mindful choice of the initial sampling data as well as of the algorithm
parameters.
Still, the optimization problem  could be solved with
a satisfactory accuracy for most of the
configurations of this extremely hard problem.

%Moreover, relative errors for the problem with moderate internal damping are   smaller compared to the problem with very small internal damping part. Therefore, the problem with smaller internal damping need to be treated more carefully and appropriate initial parameters needs to be taken into account.

\section*{Code Availability}
The code and data that was used to obtain the results in this paper is freely available. It can be downloaded from the DOI
\begin{center}
    \url{10.5281/zenodo.3634361}.
\end{center}
\section*{Acknowlegdement}
We thank Tiphaine Bonniot de Ruisselet from ENSEEIHT Toulouse (France) for performing some preliminary numerical experiments.

\bibliographystyle{plain}
\bibliography{literatur}
\end{document}

%% file: relErrviscAproblemUpdated.tikz
% This file was created by matlab2tikz.
%
%The latest updates can be retrieved from
%  http://www.mathworks.com/matlabcentral/fileexchange/22022-matlab2tikz-matlab2tikz
%where you can also make suggestions and rate matlab2tikz.
%
\begin{tikzpicture}

\begin{axis}[%
width=5.5in,
height=3.5in,
at={(0.758in,0.481in)},
scale only axis,
xmin=0,
xmax=37,
xlabel style={font=\color{white!15!black}},
xlabel={configuration number},
ymode=log,
ymin=1e-08,
ymax=10000,
yminorticks=true,
ylabel style={font=\color{white!15!black}},
ylabel={relative error in the optimal gains},
axis background/.style={fill=white},
xmajorgrids,
ymajorgrids,
legend style={at={(0.35,0.95)},legend cell align=left, align=left, draw=white!15!black}
]
\addplot [color=black, draw=none, mark size=6.0pt, mark=x, mark options={solid, black}]
  table[row sep=crcr]{%
1	1.29103993819368e-05\\
2	1.95192443983773e-06\\
3	41.9497617826847\\
4	8.74273100788555e-05\\
5	0.072632896797018\\
6	0.000278888443587038\\
7	1.19583900791391e-05\\
8	0.00586356100505142\\
9	1.23837505415734e-05\\
10	0.000136521449903601\\
11	4.82371945328888e-08\\
12	2.46969348693315e-05\\
13	4.41282089845235e-05\\
14	0.000230432562895625\\
15	0.000329844484826791\\
16	0.000365310988057191\\
17	0.0379150603135362\\
18	1.13568648366447e-08\\
19	9.61586311935528e-06\\
20	59.2453704227991\\
21	5.8688289123697e-05\\
22	4.21010952201551e-05\\
23	0.00048487262073637\\
24	0.000421501115693944\\
25	0.0304025766772724\\
26	0.0025449587909048\\
27	0.0148914727165\\
28	0.000189720948695603\\
29	0.000644610643810005\\
30	6.35319769160264e-06\\
31	3.20766776759078\\
32	5.41522242651894e-06\\
33	0.000193836401446192\\
34	0.0326176088316364\\
35	0.143508229436918\\
36	3.65191969495698e-05\\
};
\addlegendentry{\textbf{Mode i)}}

\addplot [color=green, draw=none, mark size=6.0pt, mark=o, mark options={solid, green}]
  table[row sep=crcr]{%
1	1.32660606385178e-05\\
2	1.23955490053119\\
3	35.5721073444541\\
4	6.09258972857852e-05\\
5	0.0498253840013737\\
6	0.00033408990082578\\
7	1.04174179713559e-05\\
8	0.0121759068985356\\
9	2.23186385513155e-05\\
10	8.6448443588655e-05\\
11	9.7041580840874e-07\\
12	0.00329597459882111\\
13	8.56023734316486e-06\\
14	0.000212898495396095\\
15	0.000335296469738919\\
16	0.00447837702406542\\
17	0.037917781913099\\
18	0.0786332945127878\\
19	9.60822972216916e-06\\
20	0.000155544625216205\\
21	0.000121362724757467\\
22	4.35224884075132e-05\\
23	0.000467958878155804\\
24	0.000414879918387727\\
25	0.00032864064634345\\
26	0.00290379403823924\\
27	0.000126748383068263\\
28	3.97547341191767e-05\\
29	0.000644983334388324\\
30	2.70269805616456e-07\\
31	3.18215199290047\\
32	2.33186906926442e-05\\
33	6.53200088244428e-05\\
34	4.03297619021382e-08\\
35	0.140171684087917\\
36	6.03806057484017e-05\\
};
\addlegendentry{\textbf{Mode ii)}}

\addplot [color=blue, draw=none, mark size=6.0pt, mark=asterisk, mark options={solid, blue}]
  table[row sep=crcr]{%
1	1.25973565454877e-05\\
2	1.23988622534308\\
3	0.00186920538941549\\
4	9.06384920630527e-05\\
5	0.121954611819434\\
6	0.000147687690469945\\
7	2.94529504037205e-06\\
8	0.23412579208616\\
9	4.65632119225053e-06\\
10	5.5127909001027e-05\\
11	2.22573508996443e-07\\
12	0.000182394690596737\\
13	1.38027558922298e-06\\
14	5.14780962822952e-05\\
15	0.000120240810409248\\
16	0.000305864472125182\\
17	0.0021074681948238\\
18	0.000176768966068064\\
19	1.1013276227804e-05\\
20	0.000313620639143239\\
21	9.66322151121704e-05\\
22	3.55840967454557e-05\\
23	0.000471821557927932\\
24	1.30327474557951e-05\\
25	0.00126629220267038\\
26	0.0012593392015018\\
27	0.000423489004772635\\
28	3.02699287991367e-05\\
29	0.000611528218816089\\
30	905.31836742165\\
31	1.85554514965697\\
32	5.38551073666319e-07\\
33	0.000180569579877806\\
34	8.78562476000775e-06\\
35	1.65458412710588e-05\\
36	0.00010214631896154\\
};
\addlegendentry{\textbf{Mode iii)}}

\addplot [color=red, draw=none, mark size=6.0pt, mark=triangle, mark options={solid, rotate=270, red}]
  table[row sep=crcr]{%
1	1.20769122547828e-05\\
2	1.23988590132025\\
3	0.000801884520891959\\
4	9.10278564875594e-05\\
5	0.0503941198307889\\
6	0.000152011436344271\\
7	5.37637300469853e-06\\
8	0.234157635888051\\
9	8.06778499524785e-06\\
10	1.32532843028977e-05\\
11	1.46291204002457e-06\\
12	0.000145999137525041\\
13	1.79185059065517e-06\\
14	5.41570194842911e-05\\
15	0.000120927401348729\\
16	0.000120117108138811\\
17	0.00210946241094868\\
18	0.000176794217369501\\
19	9.64542785988597e-06\\
20	0.000308494456079031\\
21	0.000119279457158168\\
22	2.58006961602607e-05\\
23	0.000471718903006832\\
24	1.2443808264069e-05\\
25	0.00114061037994397\\
26	0.00130603564127363\\
27	0.000642970277877291\\
28	3.17036571632278e-05\\
29	0.000637751002528421\\
30	872.227158689923\\
31	1.85617283787641\\
32	4.59940865511287e-06\\
33	0.000142430276823409\\
34	8.43200311605877e-06\\
35	0.000112775997084408\\
36	8.5477203227714e-05\\
};
\addlegendentry{\textbf{Mode iv)}}

\end{axis}
\end{tikzpicture}%

%% file: relErrHinfAproblemUpdated.tikz
% This file was created by matlab2tikz.
%
%The latest updates can be retrieved from
%  http://www.mathworks.com/matlabcentral/fileexchange/22022-matlab2tikz-matlab2tikz
%where you can also make suggestions and rate matlab2tikz.
%
\begin{tikzpicture}

\begin{axis}[%
width=5.5in,
height=3.5in,
at={(0.758in,0.481in)},
scale only axis,
xmin=0,
xmax=37,
xlabel style={font=\color{white!15!black}},
xlabel={configuration number},
ymode=log,
ymin=1e-12,
ymax=100,
yminorticks=true,
ylabel style={font=\color{white!15!black}},
ylabel={relative error in the function value},
axis background/.style={fill=white},
xmajorgrids,
ymajorgrids,
legend style={at={(0.568,0.746)}, anchor=south west, legend cell align=left, align=left, draw=white!15!black}
]
\addplot [color=black, draw=none, mark size=6.0pt, mark=x, mark options={solid, black}]
  table[row sep=crcr]{%
1	4.26679030802728e-10\\
2	5.22902038041852e-07\\
3	0.00238899450330464\\
4	1.87053752079043e-09\\
5	1.13152486805159e-06\\
6	3.94146458727082e-09\\
7	2.95487592342045e-06\\
8	0.00102172803014301\\
9	2.02077623615071e-08\\
10	5.42137814419351e-07\\
11	3.32555589903506e-08\\
12	1.25466531467365e-07\\
13	7.50788218899837e-09\\
14	7.64609755902946e-09\\
15	4.39287386375127e-08\\
16	5.04444661171975e-10\\
17	0.00102313504991801\\
18	8.4012608836355e-09\\
19	1.04164257136139e-07\\
20	2.00293307793438\\
21	3.8473220315545e-09\\
22	1.13746069414766e-09\\
23	2.99071294209901e-08\\
24	2.46029434024091e-06\\
25	0.000842510787819724\\
26	3.48275669409761e-09\\
27	3.1650801905465e-05\\
28	8.12156355646923e-09\\
29	1.25776189967693e-07\\
30	3.67030155627897e-09\\
31	0.179089362504839\\
32	4.99227639110191e-08\\
33	2.00392177645427e-08\\
34	0.00718945876592401\\
35	0.0187682814222844\\
% 36	3.49608613289731e-09\\
};
\addlegendentry{\textbf{Mode i)}}

\addplot [color=green, draw=none, mark size=6.0pt, mark=o, mark options={solid, green}]
  table[row sep=crcr]{%
1	4.08171314231171e-10\\
2	0.0317400766641472\\
3	0.0023741540036398\\
4	1.39193489826468e-09\\
5	1.34348811833873e-06\\
6	4.15272184257017e-09\\
7	2.94002145872802e-08\\
8	0.00101646319853524\\
9	2.7028164626333e-08\\
10	7.19596154409809e-07\\
11	1.05482190123546e-07\\
12	2.89374512328632e-06\\
13	7.06180487838699e-09\\
14	9.97016631738349e-09\\
15	4.63520968622026e-08\\
16	1.62498559572858e-09\\
17	0.00102163991801\\
18	0.000474097740937347\\
19	1.03923175587602e-07\\
20	2.2345207284096e-06\\
21	5.30834127233119e-09\\
22	9.54783475243502e-10\\
23	2.94528451206161e-08\\
24	7.28997226432666e-06\\
25	5.74745710343081e-06\\
26	9.19715533988081e-08\\
27	1.76160289822754e-07\\
28	6.11598366639535e-10\\
29	1.25665875056946e-07\\
30	3.2445187833729e-09\\
31	0.179099481392317\\
32	4.93959749064827e-08\\
33	5.74696301396503e-09\\
34	1.43389367329535e-07\\
35	0.0187606189409724\\
36	3.58619444078534e-09\\
};
\addlegendentry{\textbf{Mode ii)}}

\addplot [color=blue, draw=none, mark size=6.0pt, mark=asterisk, mark options={solid, blue}]
  table[row sep=crcr]{%
1	4.36334526432545e-10\\
2	0.0317100029560302\\
3	3.8568777385397e-09\\
4	2.63606964778758e-10\\
5	0.0471852237234405\\
6	1.84807728661943e-09\\
7	3.8383319815967e-08\\
8	0.0246162029950308\\
9	3.60419171808765e-08\\
10	1.66783010002333e-07\\
11	2.20373079547707e-08\\
12	1.63940541328535e-08\\
13	6.78806052990153e-09\\
14	1.01805977764468e-08\\
15	3.40365947865288e-09\\
16	9.26772543853118e-10\\
17	0.000277671288992885\\
18	0.000107211524019269\\
19	9.56710742111756e-08\\
20	0.000193534946885398\\
21	5.01007847746425e-09\\
22	9.22158033431881e-10\\
23	2.93708824255356e-08\\
24	2.63216897618837e-07\\
25	1.51042114934597e-06\\
26	1.45697532597636e-07\\
27	9.32022737450744e-09\\
28	3.06902044350043e-10\\
29	1.25235574154982e-07\\
30	0.49467970049205\\
31	0.0497225640763267\\
32	4.8783886062571e-08\\
33	2.82780032284245e-05\\
34	3.20237365315656e-06\\
35	6.80261116342131e-11\\
36	1.96763896391256e-09\\
};
\addlegendentry{\textbf{Mode iii)}}

\addplot [color=red, draw=none, mark size=6.0pt, mark=triangle, mark options={solid, rotate=270, red}]
  table[row sep=crcr]{%
1	4.19748912310445e-10\\
2	0.0317100085122828\\
3	3.25277279830101e-09\\
4	2.04311014170759e-10\\
5	1.34503288739331e-06\\
6	3.69576316460566e-09\\
7	3.93933433662518e-08\\
8	0.0246764728804538\\
9	3.74770441260104e-08\\
10	9.08810425215544e-06\\
11	1.26452932373698e-07\\
12	1.43498623978214e-08\\
13	6.7631014992062e-09\\
14	1.92308456010102e-08\\
15	3.64892502488666e-09\\
16	2.10048306683936e-10\\
17	0.000277757569219043\\
18	0.000107210954860455\\
19	1.04766170130934e-07\\
20	0.000193378236894667\\
21	5.20849010306589e-09\\
22	4.18758304384804e-10\\
23	2.94633118008925e-08\\
24	2.61506152999204e-07\\
25	1.48749609014927e-06\\
26	1.45230586124772e-07\\
27	1.86963318564242e-08\\
28	4.8366489929808e-10\\
29	1.25807793278262e-07\\
30	0.490936911079506\\
31	0.0496391327265011\\
32	5.0332980118992e-08\\
33	2.83391687799903e-05\\
34	3.08365584333403e-06\\
35	2.6659868943919e-08\\
36	2.77335129141314e-09\\
};
\addlegendentry{\textbf{Mode iv)}}

\end{axis}
\end{tikzpicture}%

%% file: relErrviscBproblemUpdated.tikz
% This file was created by matlab2tikz.
%
%The latest updates can be retrieved from
%  http://www.mathworks.com/matlabcentral/fileexchange/22022-matlab2tikz-matlab2tikz
%where you can also make suggestions and rate matlab2tikz.
%
\begin{tikzpicture}

\begin{axis}[%
width=5.5in,
height=3.5in,
at={(0.758in,0.481in)},
scale only axis,
xmin=0,
xmax=37,
xlabel style={font=\color{white!15!black}},
xlabel={configuration number},
ymode=log,
ymin=1e-05,
ymax=100,
yminorticks=true,
ylabel style={font=\color{white!15!black}},
ylabel={relative error in the optimal gains},
axis background/.style={fill=white},
xmajorgrids,
ymajorgrids,
legend style={at={(0.35,0.95)}, legend cell align=left, align=left, draw=white!15!black}
]
\addplot [color=black, draw=none, mark size=6.0pt, mark=x, mark options={solid, black}]
  table[row sep=crcr]{%
1	0.00421215716789138\\
2	0.0871588409413348\\
3	0.00787094286244928\\
4	0.0338615832683473\\
5	0.0171595396345236\\
6	0.00074422426219256\\
7	0.00243172173843439\\
8	0.000237273614821916\\
9	0.00378167795847253\\
10	0.0256817786284456\\
11	0.00105756361848823\\
12	0.00384356419952759\\
13	0.0265712464397314\\
14	0.00729925363030292\\
15	0.00881051749976108\\
16	0.0815898608065143\\
17	0.00622516856235869\\
18	0.744666783134425\\
19	0.016334380125797\\
20	0.0699025804231869\\
21	0.0401140830157299\\
22	0.010537192731385\\
23	0.0140756652284189\\
24	0.0142111813200278\\
25	0.0430089331956784\\
26	0.0110398006749523\\
27	0.0644773165467904\\
28	0.0102482164650838\\
29	0.0162596661012729\\
30	0.608318957199269\\
31	0.000566459406113101\\
32	0.000792379532070978\\
33	0.00890538795903989\\
34	0.110863457101741\\
35	0.0108220968183153\\
36	0.00489115566560229\\
};
\addlegendentry{\textbf{Mode i)}}

\addplot [color=green, draw=none, mark size=6.0pt, mark=o, mark options={solid, green}]
  table[row sep=crcr]{%
1	0.000152655494536892\\
2	0.00618846493219556\\
3	0.00756401857356218\\
4	0.0452483217099022\\
5	0.0235672985621215\\
6	0.000862759849837062\\
7	0.000745176255072424\\
8	0.00911058653086002\\
9	0.00395516371954105\\
10	0.0652192619400208\\
11	0.0031507122167641\\
12	0.00182669259709434\\
13	0.0113820535071777\\
14	0.0191314318777174\\
15	0.00132372540316496\\
16	0.0793043450167593\\
17	0.00568054810316308\\
18	0.00752631124085826\\
19	0.0281039804730622\\
20	0.0745250168136915\\
21	0.039533837974974\\
22	0.0127067443414981\\
23	0.0213380213742789\\
24	0.040889118268741\\
25	0.0353458264662972\\
26	0.0214508770266866\\
27	0.0028915406001074\\
28	0.267130965811519\\
29	0.0198393070899412\\
30	0.00407394373229519\\
31	0.000696479974885769\\
32	0.00620801040235987\\
33	0.0119701887710051\\
34	0.0741524297358309\\
35	0.00632687408919785\\
36	0.00489115015483198\\
};
\addlegendentry{\textbf{Mode ii)}}

\addplot [color=blue, draw=none, mark size=6.0pt, mark=asterisk, mark options={solid, blue}]
  table[row sep=crcr]{%
1	0.000171632588721124\\
2	0.0038998518950637\\
3	0.0114660730067683\\
4	0.0943151195121293\\
5	0.0084527968403925\\
6	0.00070683769573212\\
7	0.00122938752331256\\
8	0.000344775511543572\\
9	0.0104877710406782\\
10	0.0221230958589062\\
11	0.00710331886332571\\
12	0.00424474481625254\\
13	0.00281005346044563\\
14	0.0071642770216243\\
15	9.73416640559078e-05\\
16	0.0809971518359097\\
17	0.00609659344787964\\
18	0.0033073657027459\\
19	0.00708485411056349\\
20	0.06180248602858\\
21	71.1702951875286\\
22	0.0112363906889709\\
23	0.0132646665513508\\
24	0.0449837253011897\\
25	0.0248883912378331\\
26	0.029453804754898\\
27	0.000935121630354892\\
28	0.0122213140251727\\
29	0.014330615317526\\
30	0.00807799481582121\\
31	0.000825719526337628\\
32	0.000758970384889008\\
33	0.00863391730115558\\
34	0.0433448005422321\\
35	0.00504928808908864\\
36	5.83882872826424e-05\\
};
\addlegendentry{\textbf{Mode iii)}}

\addplot [color=red, draw=none, mark size=6.0pt, mark=triangle, mark options={solid, rotate=270, red}]
  table[row sep=crcr]{%
1	0.000171632638157516\\
2	0.785977331708511\\
3	0.0115422755668202\\
4	0.0944320787180362\\
5	0.00844703612822451\\
6	0.000636972168383297\\
7	0.00126581727205286\\
8	0.000126197022390285\\
9	0.0104354289126198\\
10	0.021240621905357\\
11	0.00671418784205999\\
12	0.00424474316486405\\
13	0.00295969448442504\\
14	0.00733396764199937\\
15	0.000746001306571556\\
16	0.0809971517608694\\
17	0.00609647718660767\\
18	0.00330736565261666\\
19	0.00708482051907905\\
20	0.0621650106287926\\
21	48.9577797790655\\
22	0.0137063409004455\\
23	1.71725750132343\\
24	0.0448447825179713\\
25	0.0262717000662311\\
26	0.0288429218680536\\
27	0.000935121671443757\\
28	0.0104394394476017\\
29	0.0143306156413475\\
30	0.00807799479179271\\
31	0.000825719554697081\\
32	0.000789219268329604\\
33	0.00767754539514926\\
34	0.0480301846494182\\
35	0.00537868649743158\\
36	0.000184050860952436\\
};
\addlegendentry{\textbf{Mode iv)}}

\end{axis}
\end{tikzpicture}%

%% file: relErrHinfBproblemUpdated.tikz
% This file was created by matlab2tikz.
%
%The latest updates can be retrieved from
%  http://www.mathworks.com/matlabcentral/fileexchange/22022-matlab2tikz-matlab2tikz
%where you can also make suggestions and rate matlab2tikz.
%
\begin{tikzpicture}

\begin{axis}[%
width=5.5in,
height=3.5in,
at={(0.758in,0.481in)},
scale only axis,
xmin=0,
xmax=37,
xlabel style={font=\color{white!15!black}},
xlabel={configuration number},
ymode=log,
ymin=1e-10,
ymax=0.0001,
yminorticks=true,
ylabel style={font=\color{white!15!black}},
ylabel={relative error in the function value},
axis background/.style={fill=white},
xmajorgrids,
ymajorgrids,
legend style={at={(0.9,0.3)}, legend cell align=left, align=left, draw=white!15!black}
]
\addplot [color=black, draw=none, mark size=6.0pt, mark=x, mark options={solid, black}]
  table[row sep=crcr]{%
1	3.64944709546144e-06\\
2	1.41758337282266e-06\\
3	1.92402848483707e-06\\
4	1.30619444335686e-05\\
5	1.88647376187633e-07\\
6	7.37452702121227e-08\\
7	1.03684387919863e-06\\
8	1.14650936554265e-08\\
9	1.06049017756675e-06\\
10	5.69912947046494e-06\\
11	2.81012703422609e-07\\
12	2.16239588116276e-06\\
13	1.35344295387721e-05\\
14	1.76647891169445e-06\\
15	2.29396333762943e-06\\
16	2.68667267489616e-05\\
17	1.05916679631898e-06\\
18	3.86764908766689e-08\\
19	4.00465455509682e-08\\
20	2.96756071303116e-05\\
21	1.37262457722415e-05\\
22	4.27906650119999e-07\\
23	1.58117615916893e-06\\
24	5.23174710917328e-06\\
25	5.07179226525955e-06\\
26	6.61363803465925e-06\\
27	4.78434725452865e-05\\
28	6.10011592345985e-07\\
29	3.93895575458087e-06\\
30	2.53049086799426e-06\\
31	7.70580541967985e-08\\
32	6.8106046137432e-08\\
33	1.5003739378779e-06\\
34	1.59368480942411e-06\\
35	4.86123848441426e-07\\
36	2.81433199144359e-06\\
};
\addlegendentry{\textbf{Mode i)}}

\addplot [color=green, draw=none, mark size=6.0pt, mark=o, mark options={solid, green}]
  table[row sep=crcr]{%
1	4.26015396305472e-09\\
2	1.01252799539758e-06\\
3	1.88489280580125e-06\\
4	1.6371872795729e-05\\
5	2.96188037578038e-06\\
6	6.82830825646379e-08\\
7	1.39649794528716e-07\\
8	8.47100782966193e-06\\
9	1.3524612244465e-06\\
10	2.42933481441474e-05\\
11	8.70237899474775e-07\\
12	1.4706487217681e-06\\
13	1.39500739520349e-06\\
14	5.09250166128194e-06\\
15	3.6755232969067e-08\\
16	2.60283864850591e-05\\
17	1.05459753377066e-06\\
18	4.15543284967045e-07\\
19	7.70825759332822e-06\\
20	2.92797425038063e-05\\
21	1.37104274374568e-05\\
22	8.59789219656158e-07\\
23	1.09637885651855e-06\\
24	1.21190231613998e-05\\
25	6.57571615726467e-06\\
26	9.95696553990686e-06\\
27	6.13774785232127e-08\\
28	6.24611761556031e-07\\
29	3.6313679837696e-06\\
30	5.5327976479475e-07\\
31	8.24661408940226e-08\\
32	5.07159398470083e-08\\
33	1.29871942123487e-06\\
34	6.0116523359116e-06\\
35	1.05351261615682e-06\\
36	2.81432791023847e-06\\
};
\addlegendentry{\textbf{Mode ii)}}

\addplot [color=blue, draw=none, mark size=6.0pt, mark=asterisk, mark options={solid, blue}]
  table[row sep=crcr]{%
1	6.60298490104484e-09\\
2	1.4594848802676e-06\\
3	2.13541075676023e-06\\
4	2.32101471080163e-05\\
5	1.39344204437314e-06\\
6	7.42942617157883e-08\\
7	2.21109161949855e-07\\
8	5.01634757823587e-09\\
9	2.35756017580014e-06\\
10	1.54259260573254e-06\\
11	1.29333556555782e-06\\
12	2.18309024957415e-06\\
13	3.91435162305868e-07\\
14	1.77304411480318e-06\\
15	4.41768371043732e-10\\
16	2.68387350059471e-05\\
17	1.06105899750641e-06\\
18	3.62147110192121e-07\\
19	3.68669157489313e-07\\
20	3.00679485574746e-05\\
21	1.37388572449796e-05\\
22	7.81949378178683e-07\\
23	1.5863652068227e-06\\
24	1.21989506226623e-05\\
25	6.99052128841054e-06\\
26	1.08281339447941e-05\\
27	1.53268805109201e-08\\
28	6.32279331558556e-07\\
29	3.97303156956651e-06\\
30	7.66290438297882e-07\\
31	7.99365117871909e-08\\
32	6.64087930197097e-08\\
33	1.50212325211452e-06\\
34	7.27615633751494e-06\\
35	1.03433244563974e-06\\
36	3.27010137603682e-09\\
};
\addlegendentry{\textbf{Mode iii)}}

\addplot [color=red, draw=none, mark size=6.0pt, mark=triangle, mark options={solid, rotate=270, red}]
  table[row sep=crcr]{%
1	6.61172937962083e-09\\
2	1.4591978885672e-06\\
3	2.1353681053658e-06\\
4	2.32135119422953e-05\\
5	1.39343745013308e-06\\
6	7.41034133191781e-08\\
7	2.21552085687071e-07\\
8	1.34024444636966e-09\\
9	2.35701345815404e-06\\
10	1.53644215009673e-06\\
11	1.28632168935033e-06\\
12	2.18308132634152e-06\\
13	3.91421615483268e-07\\
14	1.76954955521005e-06\\
15	2.07836777597272e-08\\
16	2.68387414436889e-05\\
17	1.06105262541989e-06\\
18	3.62132872014145e-07\\
19	3.68665599835697e-07\\
20	3.00734012237332e-05\\
21	1.37388675407452e-05\\
22	8.73737607324833e-07\\
23	1.57688927375505e-06\\
24	1.21996437013074e-05\\
25	7.04366844445791e-06\\
26	1.08121044061383e-05\\
27	1.53208553661047e-08\\
28	6.02356598337863e-07\\
29	3.97302680463374e-06\\
30	7.66299794001527e-07\\
31	7.99223607531759e-08\\
32	6.81060999711685e-08\\
33	1.48276549320756e-06\\
34	7.40888441290583e-06\\
35	1.04859970680839e-06\\
36	4.72919486809352e-09\\
};
\addlegendentry{\textbf{Mode iv)}}

\end{axis}
\end{tikzpicture}%

%% file: TomV19.bbl
\newcommand{\noopsort}[1]{} \newcommand{\printfirst}[2]{#1}
  \newcommand{\singleletter}[1]{#1} \newcommand{\switchargs}[2]{#2#1}
\begin{thebibliography}{10}

\bibitem{AliBMSV17a}
N.~Aliyev, P.~Benner, E.~Mengi, P.~Schwerdtner, and M.~Voigt.
\newblock A greedy subspace method for computing the $\mathcal{L}_\infty$-norm.
\newblock {\em PAMM Proc. Appl. Math. Mech.}, 17(1):751--752, 2017.

\bibitem{AliBMSV17}
N.~Aliyev, P.~Benner, E.~Mengi, P.~Schwerdtner, and M.~Voigt.
\newblock Large-scale computation of $\mathcal{L}_\infty$-norms by a greedy
  subspace method.
\newblock {\em SIAM J. Matrix Anal. Appl.}, 38(4):1496--1516, 2017.

\bibitem{AliBMV18}
N.~Aliyev, P.~Benner, E.~Mengi, and M.~Voigt.
\newblock A subspace framework for {$\mathcal{H}_\infty$}-norm minimization,
  2018.
\newblock In preparation.

\bibitem{Gugercin2009}
C.~Beattie and S.~Gugercin.
\newblock Interpolatory projection methods for structure-preserving model
  reduction.
\newblock {\em Systems Control Lett.}, 58(3):225--232, 2009.

\bibitem{BKTT15}
P.~Benner, P.~K{\"u}rschner, Z.~Tomljanovi\'{c}, and N.~Truhar.
\newblock {Semi-active damping optimization of vibrational systems using the
  parametric dominant pole algorithm}.
\newblock {\em ZAMM Z. Angew. Math. Mech.}, 96(4):604--619, 2016.

\bibitem{BenSV12}
P.~Benner, V.~Sima, and M.~Voigt.
\newblock $\mathcal{L}_\infty$-norm computation for continuous-time descriptor
  systems using structured matrix pencils.
\newblock {\em IEEE Trans. Automat. Control}, 57(1):233--238, 2012.

\bibitem{BennerTomljTruh10}
P.~Benner, Z.~Tomljanovi{\'c}, and N.~Truhar.
\newblock {Dimension reduction for damping optimization in linear vibrating
  systems}.
\newblock {\em ZAMM Z. Angew. Math. Mech.}, 91(3):179--191, 2011.

\bibitem{BennerTomljTruh11}
P.~Benner, Z.~Tomljanovi{\'c}, and N.~Truhar.
\newblock {Optimal damping of selected eigenfrequencies using dimension
  reduction}.
\newblock {\em Numer. Linear Algebra Appl.}, 20(1):1--17, 2013.

\bibitem{BenV13e}
P.~Benner and M.~Voigt.
\newblock A structured pseudospectral method for $\mathcal{H}_\infty$-norm
  computation of large-scale descriptor systems.
\newblock {\em Math. Control Signals Systems}, 26(2):303--338, 2014.

\bibitem{Blanchini12}
F.~Blanchini, D.~Casagrande, P.~Gardonio, and S.~Miani.
\newblock {Constant and switching gains in semi-active damping of vibrating
  structures}.
\newblock {\em Internat. J. Control}, 85(12):1886--1897, 2012.

\bibitem{BB90}
S.~Boyd and V.~Balakrishnan.
\newblock A regularity result for the singular values of a transfer matrix and
  a quadratically convergent algorithm for computing its
  {$\mathrm{L}_{\infty}$}-norm.
\newblock {\em Systems Control Lett.}, 15(1):1--7, 1990.

\bibitem{BRAB98}
K.~Brabender.
\newblock {\em {Optimale D{\"a}mpfung von linearen Schwingungssystemen}}.
\newblock Dissertation, Fachbereich Mathematik, FernUniversit{\"a}t --
  Gesamthochschule -- in Hagen, 1998.

\bibitem{BS90}
N.~A. Bruinsma and M.~Steinbuch.
\newblock A fast algorithm to compute the {$H_{\infty}$}-norm of a transfer
  function matrix.
\newblock {\em Systems Control Lett.}, 14(4):287--293, 1990.

\bibitem{BurLO05}
J.~V. Burke, A.~S. Lewis, and M.~L. Overton.
\newblock A robust gradient sampling algorithm for nonsmooth, nonconvex
  optimization.
\newblock {\em SIAM J. Optimization}, 23:751--779, 2005.

\bibitem{CurtisMitchOverton17}
F.~E. Curtis, T.~Mitchell, and M.~L. Overton.
\newblock A {BFGS-SQP} method for nonsmooth, nonconvex, constrained
  optimization and its evaluation using relative minimization profiles.
\newblock {\em Optim. Methods Softw.}, 32(1):148--181, 2017.

\bibitem{FreitLanc:99}
P.~Freitas and P.~Lancaster.
\newblock On the optimal value of the spectral abscissa for a system of linear
  oscillators.
\newblock {\em SIAM J. Matrix Anal. Appl.}, 21(1):195--208, 1999.

\bibitem{Gaw}
W.~Gawronski.
\newblock {\em {Advanced Structural Dynamics and Active Control of
  Structures}}.
\newblock Mech. Engrg. Ser. Springer, New York, USA, 1st edition, 2004.

\bibitem{GugGO13}
N.~Guglielmi, M.~G\"urb\"uzbalaban, and M.~L. Overton.
\newblock Fast approximation of the {$H_\infty$} norm via optimization over
  spectral value sets.
\newblock {\em SIAM J. Matrix Anal. Appl.}, 34(2):709--737, 2013.

\bibitem{GumHMO09}
S.~Gumussoy, D.~Henrion, M.~Millstone, and M.~L. Overton.
\newblock Multiobjective robust control with {HIFOO} 2.0.
\newblock In {\em Proc. 6th IFAC Symposium on Robust Control Design}, Haifa,
  2009.

\bibitem{Kur10}
P.~K\"urschner.
\newblock Two-sided eigenvalue algorithms for modal approximation.
\newblock Masterarbeit, Chemnitz University of Technology, Faculty of
  Mathematics, 2010.

\bibitem{Lan64}
P.~Lancaster.
\newblock On eigenvalues of matrices dependent on a parameter.
\newblock {\em Numer. Math.}, 6:377--387, 1964.

\bibitem{LewO08}
A.~S. Lewis and M.~L. Overton.
\newblock Nonsmooth optimization via {BFGS}, 2008.
\newblock Available from
  \url{https://cs.nyu.edu/overton/papers/pdffiles/bfgs_inexactLS.pdf}.

\bibitem{LewO13}
A.~S. Lewis and M.~L. Overton.
\newblock Nonsmooth optimization via quasi-{N}ewton methods.
\newblock {\em Math. Program.}, 143(1--2):135--163, 2013.

\bibitem{MS05}
V.~Mehrmann and T.~Stykel.
\newblock Balanced truncation model reduction for large-scale systems in
  descriptor form.
\newblock In P.~Benner, V.~Mehrmann, and D.~Sorensen, editors, {\em Dimension
  Reduction of Large-Scale Systems}, volume~45 of {\em Lect. Notes Comput. Sci.
  Eng.}, chapter~3, pages 89--116. Springer-Verlag, Berlin, 2005.

\bibitem{MitO16}
T.~Mitchell and M.~L. Overton.
\newblock Hybrid expansion-contraction: a robust scaleable method for
  approximating the {$H_\infty$} norm.
\newblock {\em IMA J. Numer. Anal.}, 36(3):985--1014, 2016.

\bibitem{MullerSchiehlen85}
P.~C. M{\"u}ller and W.~Schiehlen.
\newblock {\em {Linear Vibrations}}, volume~7 of {\em Mechanics: Dynamical
  Systems}.
\newblock Springer Netherlands, 1985.

\bibitem{RomM08}
J.~Rommes and N.~Martins.
\newblock Computing transfer function dominant poles of large-scale
  second-order dynamical systems.
\newblock {\em SIAM J. Sci. Comput.}, 30(4):2137--2157, 2008.

\bibitem{SteS90}
G.~W. Stewart and J.-G. Sun.
\newblock {\em Matrix Perturbation Theory}.
\newblock Computer Science and Scientific Computing. Academic Press, 1990.

\bibitem{TomljGugerBeatt17}
Z.~Tomljanovi{\'c}, C.~Beattie, and S.~Gugercin.
\newblock Damping optimization of parameter dependent mechanical systems by
  rational interpolation.
\newblock {\em Adv. Comput. Math.}, 44(6):1797--1820, 2018.

\bibitem{TrTomPuv16}
N.~Truhar, Z.~Tomljanovi\'{c}, and M.~Puva\v{c}a.
\newblock An efficient approximation for optimal damping in mechanical systems.
\newblock {\em Int. J. Numer. Anal. Model.}, 14(2):201--217, 2017.

\bibitem{TRUTOMVES2014}
N.~Truhar, Z.~Tomljanovi\'{c}, and K.~Veseli\'{c}.
\newblock Damping optimization in mechanical systems with external force.
\newblock {\em Appl. Math. Comput.}, 250:270--279, 2015.

\bibitem{TRUHVES09}
N.~Truhar and K.~Veseli{\'c}.
\newblock {An efficient method for estimating the optimal dampers' viscosity
  for linear vibrating systems using {Lyapunov} equation}.
\newblock {\em SIAM J. Matrix Anal. Appl.}, 31(1):18--39, 2009.

\bibitem{VES2011}
K.~Veseli{\'c}.
\newblock {\em {Damped Oscillations of Linear Systems}}, volume 2023 of {\em
  Lecture Notes in Math.}
\newblock Springer, Berlin, Heidelberg, 2011.

\bibitem{Zhou96}
K.~Zhou, J.~C. Doyle, and K.~Glover.
\newblock {\em {Robust and Optimal Control}}.
\newblock Prentice Hall, Upper Saddle River, NJ, 1996.

\end{thebibliography}
